\documentclass{article}
\usepackage{amssymb}
\title{\bf
Multiple elliptic hypergeometric series 
\\ 
{\large --- An approach from the Cauchy determinant ---}
}
\author{Yasushi KAJIHARA\footnote{
Department of Mathematics, Graduate School of Science, Osaka University}
\ and Masatoshi NOUMI\footnote{
Department of Mathematics, Graduate School of Science and Technology, 
Kobe University}}
\date{}
\newtheorem{thm}{Theorem}[section]
\newtheorem{prop}[thm]{Proposition}
\newtheorem{lem}[thm]{Lemma}
\newtheorem{rem}[thm]{Remark}
\newcommand{\dfrac}[2]{{\displaystyle\frac{#1}{#2}}}
\newcommand{\dsum}[2]{{\displaystyle\sum_{#1}^{#2}}}
\newcommand{\dprod}[2]{{\displaystyle\prod_{#1}^{#2}}}
\newcommand{\eqref}[1]{$(${\rm\ref{#1}}$)$}

\newcommand{\comment}[1]{}
\newenvironment{eqs}[1]
{\vspace{-2pt}\begin{equation}\arraycolsep=2pt\begin{array}{#1}}
{\vspace{-2pt}\end{array}\end{equation}\hspace{-4pt}}

\newcommand{\MP}[4]{
\Phi^{#1}_{#2}\Big(
\arraycolsep=1pt
\renewcommand{\arraystretch}{0.8}
\begin{array}{cccccccc}#3\end{array}\Big|
\begin{array}{cccccccc}#4\end{array}
\renewcommand{\arraystretch}{1.0}
\Big)
}
\newcommand{\Mphi}[4]{
\phi^{#1}_{#2}\Big(
\arraycolsep=1pt
\renewcommand{\arraystretch}{0.8}
\begin{array}{cccccccc}#3\end{array}\Big|
\begin{array}{cccccccc}#4\end{array}
\renewcommand{\arraystretch}{1.0}
\Big)
}
\newcommand{\Mphiu}[4]{
\phi^{#1}\Big(
\arraycolsep=1pt
\renewcommand{\arraystretch}{0.8}
\begin{array}{cccccccc}#2\end{array}\Big|
\begin{array}{cccccccc}#3\end{array}; {#4}
\renewcommand{\arraystretch}{1.0}
\Big)
}
\newcommand{\ME}[5]{
E^{#1}\Big(
\renewcommand{\arraystretch}{0.8}
\begin{array}{cccccccc}#2\end{array}
\renewcommand{\arraystretch}{1.0}
\Big|\,{#3};{#4};{#5}\Big)
}
\newcommand{\MW}[5]{
W^{#1}\Big(
\renewcommand{\arraystretch}{0.8}
\begin{array}{cccccccc}#2\end{array}
\renewcommand{\arraystretch}{1.0}
\Big|\,{#3};{#4};{#5}\Big)
}
\newcommand{\mat}[1]{
\renewcommand{\arraystretch}{0.8}
\begin{array}{cccccccccc}#1\end{array}
\renewcommand{\arraystretch}{1.0}
}

\begin{document}
\maketitle
\centerline{\sl Dedicated to Professor Tom Koornwinder 
on his sixtieth birthday}

\begin{abstract}
A multiple generalization of elliptic hypergeometric 
series is studied through the Cauchy determinant 
for the Weierstrass sigma function. 
In particular, a duality transformation for multiple 
hypergeometric series is proposed. 
As an application, two types of Bailey transformations 
for very well-poised multiple elliptic hypergeometric 
series are derived. 
\end{abstract}

\centerline{
}

\section*{Introduction}

In this paper we investigate a multiple generalization of 
elliptic hypergeometric series, and propose a 
{\em duality transformation} for 
multiple hypergeometric series. 
Our duality transformation is obtained from  
an identity arising from the Cauchy determinant 
formula for the Weierstrass sigma function, 
by means of specialization of a particular form. 
This formula for the multiple elliptic case, as well as the idea of proof, 
is a variant of the one previously studied by one of the authors 
\cite{K, K2} 
for multiple basic hypergeometric series. 
As an application of the duality transformation, we 
prove two types of Bailey transformation formulas 
for very well-poised multiple elliptic hypergeometric series; 
one of them is 
the one proved recently by Rosengren \cite{R2}, but 
the other seems to be new. 

We remark that the Bailey transformation formula for 
very well-poised elliptic 
hypergeometric series ${}_{10}E_{9}$ 
was discovered  
by Frenkel-Turaev \cite{FT} in the context of 
elliptic $6$\,-$j$ symbols.
It plays a crucial role in the theory of 
Spiridonov-Zhedanov \cite{SZ} on biorthogonal rational functions 
on elliptic grids. 
Also, in the geometric approach to the elliptic Painlev\'e 
equation \cite{KMNOY}, this transformation formula 
is explained as the standard 
Cremona transformation applied to its hypergeometric solutions.

The Bailey transformation formula \cite{B1, B2} has been 
extended 
by the method of Bailey chains 
into various directions
including multiple generalizations
(see \cite{AAR}, \cite{MN} for instance). 
We expect that our approach, based on the Cauchy determinant 
and the duality transformation, would provide a unified and 
transparent way of looking at the symmetry structure 
of Bailey transformations for the hierarchy of elliptic, 
basic and ordinary hypergeometric series. 
Multiple elliptic hypergeometric series are also generalized 
by means of root systems \cite{R2}. 
In terms of root systems, our discussion is restricted to 
the cases of $A$ type. It would be an intriguing problem to 
extend the method of Cauchy determinants to other root 
systems. 

\par\smallskip\noindent
{\small{\bf Notes.}
After completing this work, the authors were informed by 
H.~Rosengren that he obtained independently the duality 
transformation formula (Theorem 2.1). 
}

\section{Starting from the Cauchy formula}
\setcounter{equation}{0}
Throughout this paper, by the symbol 
$\,[x]$ we denote a nonzero holomorphic odd function on 
$\mathbb{C}$ in the variable $x$, 
satisfying the Riemann relation: 
\begin{eqs}{ll}\smallskip
(1) \quad&\,[-x]=-\,[x],
\\ \smallskip
(2) \quad& \,[x+y]\,[x-y]\,[u+v]\,[u-v]\\
&=\,[x+u]\,[x-u]\,[y+v]\,[y-v]
-\,[x+v]\,[x-v]\,[y+u]\,[y-u].
\end{eqs}
There are three classes of such functions; 
elliptic, trigonometric and rational. 
It is well known \cite{WW} that all $[x]$ satisfying 
the conditions above are obtained from the 
three functions 
$\sigma(x;\omega_1,\omega_2)$, 
$\sin(\pi x)$ and $x$ by the transformation 
of the form $e^{a x^2+b} [c x]$ ($a,b,c\in\mathbb{C}$);
$\sigma(x;\omega_1,\omega_2)$ stands for the 
Weierstrass sigma function 
with respect to the periods $(\omega_1,\omega_2)$.
The following determinant formula for $[x]$ plays 
a crucial role in our arguments. 
\begin{thm}[Cauchy determinant formula]\label{thm:Cauchy}
For two sets of variables 
$z=(z_1,\ldots,z_M)$ and $w=(w_1,\ldots,w_M)$, 
we have 
\begin{equation}\label{eq:Fay}
\det\left(
\dfrac{\,[\lambda+z_i+w_j]}{\,[\lambda]\,[z_i+w_j]}
\right)_{i,j=1}^M
=
\dfrac{\,\big[\lambda+\dsum{i=1}{M}(z_i+w_i)\big]\ 
\dprod{1\le i<j\le M}{}\,[z_i-z_j]\,[w_i-w_j]}
{\,[\lambda]\,\dprod{1\le i,j\le M}{}\,[z_i+w_j]},
\end{equation}
where $\lambda$ is a complex parameter. 
\end{thm}
This determinant formula can be regarded as the special 
case of Fay's trisecant formula \cite{Fay}, concerning 
the cases of elliptic curves and their degenerations. 
We remark that the elliptic case can already be found 
in Frobenius \cite{Fr}.  
Also, from the rational case where $[x]=x$, by taking the 
limit  $\lambda\to\infty$ we obtain the usual Cauchy 
determinant formula. 
Note that the case of $M=2$ is precisely the Riemann 
relation; the general formula 
\eqref{eq:Fay} can be proved by the induction on $M$, 
thanks to the formula 
\begin{equation}
\det\big(a_{ij}\big)_{i,j=1}^{M}
=a_{MM}^{2-M}\det\big(a_{ij}a_{MM}-a_{iM} a_{Mj}\big)_{i,j=1}^{M-1}
\end{equation}
for any $M\times M$ matrix $\big(a_{ij}\big)_{i,j}^M$ with $a_{MM}\ne 0$. 
Regarding the Cauchy determinant \eqref{eq:Fay} as 
a function in $z=(z_1,\ldots,z_M)$ and $w=(w_1,\ldots,w_M)$, 
we set 
\begin{equation}
D(z|w)=\det\left(
\dfrac{\,[\lambda+z_i+w_j]}{\,[\lambda]\,[z_i+w_j]}
\right)_{i,j=1}^M. 
\end{equation}
Note that $D(z|w)$ is alternating 
under the permutation of variables 
$z_i$ (resp. $w_j$), and symmetric 
with respect to exchanging $z$ and $w$. 

\par\medskip
We now fix a generic constant $\delta\in\mathbb{C}$ 
such that $[k\delta]\ne 0$ for any nonzero $k\in\mathbb{Z}$; 
in the case of the sigma function 
$[x]=\sigma(x;\omega_1,\omega_2)$, 
this condition is equivalent to 
$\delta\not\in\mathbb{Q}\omega_1+\mathbb{Q}\omega_2$. 
We denote by $T_{z_i}$ the shift operator in the variable $z_i$ 
by $\delta$: 
\begin{equation}
T_{z_i}f(z_1,\ldots,z_M)=f(z_1,\ldots,z_i+\delta,\ldots,z_M). 
\end{equation}
Given a set of indices $K=\{i_1,\ldots,i_r\}\subset\{1,\ldots,M\}$, 
we also define 
\begin{equation}
T_{z,K}=\dprod{i\in K}{}T_{z_i}=T_{z_{i_1}}\,\cdots\,T_{z_{i_r}}. 
\end{equation}
In this notation, the operator 
$E(T_z;u)=(1+u\,T_{z_1})\,\cdots\,(1+u\,T_{z_M})$, 
containing a formal parameter $u$, is expanded into the sum  
\begin{equation}\label{eq:Eexp}
E(T_z;u)=(1+u\,T_{z_1})\,\cdots\,(1+u\,T_{z_M})
=\dsum{K\subset\{1,\ldots,M\}}\ u^{|K|}\, T_{z,K}
\end{equation}
over all subsets $K$ of the index set $\{1,\ldots,M\}$.
By applying the operator $E(T_z;u)$ to the Cauchy determinant $D(z|w)$ 
we obtain 
\begin{eqs}{rl}
E(T_z;u)\,D(z|w)&=\det\left((1+u\,T_{z_i})
\dfrac{\,[\lambda+z_i+w_j]}{\,[\lambda]\,[z_i+w_j]}\right)_{i,j=1}^M\\
&=\det\left(
\dfrac{\,[\lambda+z_i+w_j]}{\,[\lambda]\,[z_i+w_j]}
+ u \,
\dfrac{\,[\lambda+z_i+w_j+\delta]}{\,[\lambda]\,[z_i+w_j+\delta]}
\right)_{i,j=1}^M.
\end{eqs}
Since this formula is symmetric with respect 
to exchanging $z$ and $w$, we have 
\begin{equation}
E(T_z;u)\,D(z|w)=E(T_w;u)\,D(w|z). 
\end{equation}
As we will see below, 
this formula, combined with the expansion \eqref{eq:Eexp}, 
gives rise to remarkable identities concerning the function 
$[x]$.  

By using the notation 
\begin{equation}
\Delta(z)=\dprod{1\le i<j\le M}{}\,[z_i-z_j],\quad |z|=z_1+\cdots+z_M,
\end{equation}
let us write the Cauchy formula in the form
\begin{equation}
D(z|w)=
\dfrac{\,[\,\lambda+|z|+|w|\,]\ 
\Delta(z)\,\Delta(w)}
{\,[\lambda]\,\dprod{1\le i,j\le M}{}\,[z_i+w_j]}. 
\end{equation}
By applying the shift operator $T_{z,K}$ to the 
difference product $\Delta(z)$, we obtain 
\begin{equation}
\dfrac{T_{z,K}\Delta(z)}{\Delta(z)}
=\dprod{i\in K; j\not\in K}\,\dfrac{[z_i-z_j+\delta]}{[z_i-z_j]},
\end{equation}
and hence 
\begin{eqs}{ll}\smallskip
\dfrac{T_{z,K} D(z|w)}{D(z|w)}
=&\dfrac{[\,\lambda+|z|+|w|+|K|\delta\,]}
{[\,\lambda+|z|+|w|\,]}
\\ 
&\cdot\ 
\dprod{i\in K;\,j\not\in K}{}\,
\dfrac{[z_i-z_j+\delta]}{[z_i-z_j]}
\dprod{i\in K;\,1\le k\le M}{}\,
\dfrac{[z_i+w_k]}{[z_i+w_k+\delta]}. 
\end{eqs}
This computation implies the following expansion of  $E(T_z;u)D(z|w)$: 
\begin{eqs}{rl}\medskip
&\dfrac{E(T_z;u)D(z|w)}{D(z|w)}
=
\dsum{K\subset\{1,\ldots,M\}}{}\,
u^{|K|}\dfrac{T_{z,K} D(z|w)}{D(z|w)}\\
&
\qquad=
\dsum{K\subset\{1,\ldots,M\}}{}\,u^{|K|}\,
\dfrac{[\,\lambda+|z|+|w|+|K|\delta\,]}
{[\,\lambda+|z|+|w|\,]}\\
&\qquad\qquad\quad\cdot
\dprod{i\in K;\,j\not\in K}{}\,
\dfrac{[z_i-z_j+\delta]}{[z_i-z_j]}
\dprod{i\in K;\,1\le k\le M}{}\,
\dfrac{[z_i+w_k]}{[z_i+w_k+\delta]}. 
\end{eqs}
In view of this formula, we denote the 
right-hand side by $F(z|w;u)$: 
\begin{eqs}{ll}\label{eq:defF}
F(z|w;u)
&=\dsum{K\subset\{1,\ldots,M\}}{}\,u^{|K|}\,
\dfrac{[\,\lambda+|z|+|w|+|K|\delta\,]}
{[\,\lambda+|z|+|w|\,]}\ 
\\
&\qquad\quad\cdot
\dprod{i\in K;\,j\not\in K}{}\,
\dfrac{[z_i-z_j+\delta]}{[z_i-z_j]}
\dprod{i\in K;\,1\le k\le M}{}\,
\dfrac{[z_i+w_k]}{[z_i+w_k+\delta]}. 
\end{eqs}
{}From the symmetry $E(T_z;u)D(z;w)=E(T_w;u)D(w;z)$,
we conclude that the function 
$F(z|w;u)$ defined as above is symmetric 
with respect to exchanging $z$ and $w$.

\begin{thm}\label{thm:FF} 
For two sets of variables $z=(z_1,\ldots,z_M)$ and 
$w=(w_1,\ldots,w_M)$,
define the function $F(z|w;u)$ by \eqref{eq:defF}. 
Then $F(z|w;u)$ is symmetric with respect to $z$ and $w$, namely,  
$F(z|w;u)=F(w|z;u)$.
\end{thm}
We remark that $F(z|w;u)$ is characterized by the identity
\begin{eqs}{ll}\medskip
E(T_z;u)D(z|w)
&=
\det\left(
\dfrac{\,[\lambda+z_i+w_j]}{\,[\lambda]\,[z_i+w_j]}
+ u \,
\dfrac{\,[\lambda+z_i+w_j+\delta]}{\,[\lambda]\,[z_i+w_j+\delta]}
\right)_{i,j=1}^M
\\
&=\ F(z|w;u)\ \dfrac{[\lambda+|z|+|w|]\ \Delta(z)\,\Delta(w)}
{[\lambda]\dprod{1\le i,j\le M}{}[z_i+w_j]}.
\end{eqs}

Explicitly, Theorem \ref{thm:FF} implies the identity 
\begin{eqs}{ll}\label{eq:FF}\smallskip
&\dsum{K\subset\{1,\ldots,M\}}{}\,u^{|K|}\,
\dfrac{[\,\lambda+|z|+|w|+|K|\delta\,]}
{[\,\lambda+|z|+|w|\,]}
\\ \medskip
&\quad\cdot\ \dprod{i\in K;\,j\not\in K}{}\,
\dfrac{[z_i-z_j+\delta]}{[z_i-z_j]}
\dprod{i\in K;\,1\le k\le M}{}\,
\dfrac{[z_i+w_k]}{[z_i+w_k+\delta]}\\
=&
\dsum{L\subset\{1,\ldots,M\}}{}\,u^{|L|}\,
\dfrac{[\,\lambda+|z|+|w|+|L|\delta\,]}
{[\,\lambda+|z|+|w|\,]}
\\ \smallskip
&\quad\cdot\ \dprod{k\in L;\,l\not\in L}{}\,
\dfrac{[w_k-w_l+\delta]}{[w_k-w_l]}
\dprod{1\le i\le M;\, k\in L}{}\,
\dfrac{[z_i+w_k]}{[z_i+w_k+\delta]}. 
\end{eqs}
By taking the coefficients of $u^d$
$(d=0,1,\ldots,M)$, 
we obtain 

\begin{thm}\label{thm:FFd}
Given two sets of variables 
$z=(z_1,\ldots,z_M)$ and $w=(w_1,\ldots,w_M)$, 
the following identity holds for each $d=0,1,\ldots,M$\,$:$ 
\begin{eqs}{ll}\smallskip
&\dsum{|K|=d}{}\ 
\dprod{i\in K;\,j\not\in K}{}\,
\dfrac{[z_i-z_j+\delta]}{[z_i-z_j]}
\dprod{i\in K;\,1\le k\le M}{}\,
\dfrac{[z_i+w_k]}{[z_i+w_k+\delta]}\\
=&
\dsum{|L|=d}{}\ 
\dprod{k\in L;\,l\not\in L}{}\,
\dfrac{[w_k-w_l+\delta]}{[w_k-w_l]}
\dprod{1\le i\le M;\, k\in L}{}\,
\dfrac{[w_k+z_i]}{[w_k+z_i+\delta]},
\end{eqs}
where $K$ and $L$ run over all $d$-subsets of 
$\{1,\ldots,M\}$.
\end{thm}

\section{Passage to multiple hypergeometric series}
\setcounter{equation}{0}
{}From $F(z|w;u)$ introduced above,  
one can generate a class of multiple hypergeometric series, 
denoted by $\Phi^{m,n}_N$ below, 
by a particular specialization of the variables $z$ and $w$. 
The symmetry property $F(z|w;u)=F(w|z;u)$ is then 
translated into a transformation formula between 
$\Phi^{m,n}_N$ and $\Phi^{n,m}_N$
(Theorem \ref{thm:PPN}).
\par\medskip

Taking $m$ variables $x=(x_1,\ldots,x_m)$ 
and a multi-index $\alpha=(\alpha_1,\ldots,\alpha_m)\in\mathbb{N}^m$
such that $|\alpha|=\alpha_1+\cdots+\alpha_m=M$,
we specialize the variables $z=(z_1,\ldots,z_M)$ as follows:
\begin{eqs}{ll}\label{eq:zsp}
(z_1,\ldots,z_M)
&=(-x_1,-x_1-\delta,\ldots,-x_1-(\alpha_1-1)\delta;
-x_2,-x_2-\delta,\ldots\\
&\qquad\ldots;-x_m,-x_m-\delta,\ldots,-x_m-(\alpha_m-1)\delta). 
\end{eqs}
Namely, we divide the index set 
$\{1,\ldots,M\}$ into $m$ blocks of size $\alpha_1,\ldots,\alpha_m$ 
according to the multi-index $\alpha$,
and then specialize each block of $z$ variables 
by the sequence $(-x_r,-x_r-\delta,\ldots,-x_r-(\alpha_r-1)\delta)$
for $r=1,\ldots,m$. 
Similarly, taking $n$ variables $y=(y_1,\ldots,y_n)$ and 
a multi-index $\beta\in\mathbb{N}^n$ with $|\beta|=M$, 
we specialize the $w$ variables as follows: 
\begin{eqs}{ll}\label{eq:wsp}
(w_1,\ldots,w_M)
&=(-y_1,-y_1-\delta,\ldots,-y_1-(\beta_1-1)\delta;
-y_2,-y_2-\delta,\ldots\\
&\qquad\ldots; -y_n,-y_n-\delta,\ldots,-y_n-(\beta_n-1)\delta).
\end{eqs}
In the following, we denote these specializations simply by 
$z\to - (x)_\alpha$, $w\to -(y)_\beta$, respectively. 
We first look at the factor 
\begin{equation}\label{eq:fct}
\dprod{i\in K;\,j\not\in K}{}\,
\dfrac{[z_i-z_j+\delta]}{[z_i-z_j]}
\end{equation}
in the expansion \eqref{eq:defF} of $F(z|w;u)$. 
When we perform the specialization $z\to-(x)_\alpha$ 
as in \eqref{eq:zsp}, 
this factor becomes zero 
if $z_j=z_i+\delta$ for some $(i,j)$ such that 
$i\in K$, $j\not\in K$. 
Replace the index set $\{1,\ldots,M\}$ for the $z$ variables 
by 
\begin{equation}
\{(r,\mu) \ | \ 1\le r \le m \ ;\ 0\le \mu <\alpha_r\}
\end{equation}
according to the $m$ blocks of $x$ variables, 
so that 
$z_{(r,\mu)}=-x_r-\mu\delta.$
Then by the specialization $z\to-(x)_\alpha$ of \eqref{eq:zsp}, 
the factor \eqref{eq:fct} becomes zero if 
there exists $(r,\mu)$ such that 
$(r,\mu)\not\in K$, $(r,\mu+1)\in K$.
This implies that, in the summation \eqref{eq:defF}, 
nontrivial terms arise only from such $K$ that all the 
elements of $K$ are packed towards  
the left in each block. 
Such subsets $K$ can be parametrized by 
the multi-indices $\gamma\in\mathbb{N}^m$ 
such that $\gamma\le\alpha$
(i.e., $\gamma_r\le\alpha_r$ for each $r=1,\ldots,m$ )
as 
\begin{equation}
K=\{(r,\mu)\ |\ 1\le r\le m\ ;\ 0\le \mu<\gamma_r\}. 
\end{equation}
With this parametrization of $K$, we compute 
\begin{eqs}{ll}\smallskip
\dprod{i\in K;\,j\not\in K}{}\,
\dfrac{[z_i-z_j+\delta]}{[z_i-z_j]}\Big|_{z\to-(x)_\alpha}\\ \smallskip
=\dprod{1\le r,s\le m}{}\ 
\dprod{0\le \mu<\gamma_r}{}\ 
\dprod{\gamma_s\le\nu <\alpha_s}{}\ 
\dfrac{[-x_r+x_s+(-\mu+\nu+1)\delta]}
{[-x_r+x_s+(-\mu+\nu)\delta]}\\ \smallskip
=\dprod{1\le r,s\le m}{}\ 
\dprod{0\le \mu<\gamma_r}{}\ 
\dfrac{[-x_r+x_s+(-\mu+\alpha_s)\delta]}
{[-x_r+x_s+(-\mu+\gamma_s)\delta]}\\
=\dprod{1\le r,s\le m}{}\ 
\dfrac{[x_r-x_s-\alpha_s\delta]_{\gamma_r}}
{[x_r-x_s-\gamma_s\delta]_{\gamma_r}},
\end{eqs}
where we have used the notation of shifted factorials
\begin{equation}
[x]_k=[x][x+\delta]\cdots [x+(k-1)\delta]\qquad(k=0,1,2,\ldots). 
\end{equation}
By the specialization $z\to -(x)_\alpha$ 
and $w\to -(y)_\beta$,
the remaining factor leads to  
\begin{eqs}{ll}\smallskip
\dprod{i\in K;\,1\le k\le M}{}\,
\dfrac{[z_i+w_k]}{[z_i+w_k+\delta]}
\Big|_{z\to-(x)_\alpha, w\to-(y)_\beta}\\ \smallskip
=\dprod{r=1}{m}\ \dprod{s=1}{n}\ 
\dprod{0\le \mu <\gamma_r}{}\ \dprod{0\le \nu<\beta_s}{}
\dfrac{[-x_r-y_s+(-\mu-\nu)\delta]}
{[-x_r-y_s+(-\mu-\nu+1)\delta]}\\\smallskip
=\dprod{r=1}{m}\ \dprod{s=1}{n}\  
\dprod{0\le \mu<\gamma_r}{}
\dfrac{[-x_r-y_s+(-\mu-\beta_s+1)\delta]}
{[-x_r-y_s+(-\mu+1)\delta]}\\\smallskip
=\dprod{r=1}{m}\ \dprod{s=1}{n}\ 
\dfrac{[x_r+y_s+(\beta_s-1)\delta]_{\gamma_r}}
{[x_r+y_s-\delta]_{\gamma_r}}. 
\end{eqs}
After all, we have 
\begin{eqs}{ll}\label{eq:Fxy}\smallskip
&F(z|w;u)\Big|_{z\to -(x)_\alpha,w\to-(y)_\beta} 
\\ \smallskip
&=\dsum{\gamma}{}\  u^{|\gamma|} \ 
\dfrac{[\lambda-\alpha x-\beta y-
{\alpha\choose 2}\delta-{\beta\choose 2}\delta
+|\gamma|\delta]}
{[\lambda-\alpha x-\beta y-
{\alpha\choose 2}\delta-{\beta\choose 2}\delta]}
\\
&\quad\cdot
\dprod{1\le r,s\le n}{}\ 
\dfrac{[x_r-x_s-\alpha_s\delta]_{\gamma_r}}
{[x_r-x_s-\gamma_s\delta]_{\gamma_r}}\ 
\dprod{1\le r\le m; 1\le s\le n}{}\ 
\dfrac{[x_r+y_s+(\beta_s-1)\delta]_{\gamma_r}}
{[x_r+y_s-\delta]_{\gamma_r}},
\end{eqs}
summed over all $\gamma\in\mathbb{N}^m$ 
with $\gamma\le\alpha$. 
Here we have used the abbreviation 
\begin{equation}
\alpha x=\alpha_1x_1+\cdots+\alpha_m x_m,
\quad
{\alpha\choose 2}={\alpha_1\choose 2}+
\cdots+{\alpha_m\choose 2}
\end{equation}
for $\alpha=(\alpha_1,\ldots,\alpha_m)$. 
In the summation \eqref{eq:Fxy}, 
notice that 
\begin{eqs}{ll}
\dprod{1\le r,s\le m}{}\ 
\dfrac{[x_r-x_s-\alpha_s\delta]_{\gamma_r}}
{[x_r-x_s-\gamma_s\delta]_{\gamma_r}}
&=
\dprod{r=1}{m}\ 
\dfrac{[-\alpha_r\delta]_{\gamma_r}}
{[-\gamma_r\delta]_{\gamma_r}}
\dprod{r\ne s}{}\ 
\dfrac{[x_r-x_s-\alpha_s\delta]_{\gamma_r}}
{[x_r-x_s-\gamma_s\delta]_{\gamma_r}}. 
\end{eqs}
Since this factor becomes zero for 
$\gamma\in\mathbb{N}^m$ such that 
$\gamma\not\le\alpha$, 
we can consider the summation 
\eqref{eq:Fxy} as being taken over all 
$\gamma\in\mathbb{N}^m$.
For convenience, we rewrite this formula further 
by using the equality
\begin{eqs}{ll}\smallskip
&\dprod{i,j}{}\ \dfrac{1}{[x_i-x_j-\gamma_j\delta]_{\gamma_i}}
\\ \smallskip
=&(-1)^{|\gamma|}\,\dprod{i<j}{}\ \dfrac{
[x_i-x_j+(\gamma_i-\gamma_j)\delta]}{[x_i-x_j]}\,
\dprod{i,j}{}\ \dfrac{1}{[x_i-x_j+\delta]_{\gamma_i}}
\\ 
=&(-1)^{|\gamma|}\,
\dfrac{\Delta(x+\gamma\delta)}{\Delta(x)}\,
\dprod{i,j}{}\ \dfrac{1}{[x_i-x_j+\delta]_{\gamma_i}},
\end{eqs}
so that
\begin{eqs}{ll}\label{eq:Fxy3}\smallskip
&F(z|w;u)\Big|_{z\to-(x)_\alpha,w\to-(y)_\beta}
\\ \smallskip
&=\dsum{\gamma\in\mathbb{N}^m}{}\  (-u)^{|\gamma|} \ 
\dfrac{[\lambda-\alpha x-\beta y-
{\alpha\choose 2}\delta-{\beta\choose 2}\delta
+|\gamma|\delta]}
{[\lambda-\alpha x-\beta y-
{\alpha\choose 2}\delta-{\beta\choose 2}\delta]}\,
\dfrac{\Delta(x+\gamma\delta)}{\Delta(x)}\,
\\
&\qquad\cdot
\dprod{1\le i,j\le m}{}\,
\dfrac{[x_i-x_j-\alpha_j\delta]_{\gamma_i}}
{[x_i-x_j+\delta]_{\gamma_i}}\ 
\dprod{1\le i\le m;\,1\le k\le n}{}\,
\dfrac{[x_i+y_k+(\beta_k-1)\delta]_{\gamma_i}}
{[x_i+y_k-\delta]_{\gamma_i}}. 
\end{eqs}

\par\medskip
We now specialize the equality 
$F(z|w;u)=F(w|z;u)$ by $z\to-(x)_\alpha, w\to-(y)_\beta$,
and take the coefficients of $u^d$ 
($d=0,1,\ldots,M$).  
Then, for any multi-indices 
$\alpha\in\mathbb{N}^m$, $\beta\in\mathbb{N}^n$ 
such that $|\alpha|=|\beta|$, 
we obtain the identity 
\begin{eqs}{ll}\smallskip
\dsum{\mu\in\mathbb{N}^m; |\mu|=d}{}\ 
\dfrac{\Delta(x+\mu\delta)}{\Delta(x)}
\dprod{i,j}{}\ 
\dfrac{[x_i-x_j-\alpha_j\delta]_{\mu_i}}
{[x_i-x_j+\delta]_{\mu_i}}\ 
\dprod{i,k}{}
\dfrac{[x_i+y_k+(\beta_k-1)\delta]_{\mu_i}}
{[x_i+y_k-\delta]_{\mu_i}}\\
=\dsum{\nu\in\mathbb{N}^n; |\nu|=d}{}\ 
\dfrac{\Delta(y+\nu\delta)}{\Delta(y)}
\dprod{k,l}{}\ 
\dfrac{[y_k-y_l-\beta_l\delta]_{\nu_k}}
{[y_k-y_l+\delta]_{\nu_k}}\ 
\dprod{k,i}{}
\dfrac{[y_k+x_i+(\alpha_i-1)\delta]_{\nu_k}}
{[y_k+x_i-\delta]_{\nu_k}}
\end{eqs}
for each $d=0,1,\ldots,M$,
where $i,j\in\{1,\ldots,m\}$ and $k,l\in\{1,\ldots,n\}$.
Replacing $y_k$ by $y_k+\delta$ ($k=1,\ldots,n$) 
in this identity, we obtain 
\begin{thm}\label{thm:PPd}
Take two sets of variables 
$x=(x_1,\ldots,x_m)$, $y=(y_1,\ldots,y_n)$
and two multi-indices 
$\alpha=(\alpha_1,\ldots,\alpha_m)\in \mathbb{N}^m$, 
$\beta=(\beta_1,\ldots,\beta_n)\in \mathbb{N}^n$ 
such that $|\alpha|=|\beta|$, i.e.,
\begin{eqs}{l}
\alpha_1+\cdots+\alpha_m=\beta_1+\cdots+\beta_n.
\end{eqs} 
Then the following identity holds for each 
$d=0,1,2,\ldots$\,$:$
\begin{eqs}{ll}\label{eq:xy}
\medskip
&\dsum{\mu\in\mathbb{N}^m;\,|\mu|=d}{}\ 
\dfrac{\Delta(x+\mu\delta)}{\Delta(x)}\, 
\dprod{i,j}{}\ 
\dfrac{[x_i-x_j-\alpha_j\delta]_{\mu_i}}
{[x_i-x_j+\delta]_{\mu_i}}\ 
\dprod{i,k}{}\ 
\dfrac{[x_i+y_k+\beta_k\delta]_{\mu_i}}
{[x_i+y_k]_{\mu_i}}\\
=&\dsum{\nu\in\mathbb{N}^n;\,|\nu|=d}{}\  
\dfrac{\Delta(y+\nu\delta)}{\Delta(y)}\, 
\dprod{k,l}{}\ 
\dfrac{[y_k-y_l-\beta_l\delta]_{\nu_k}}
{[y_k-y_l+\delta]_{\nu_k}}\ 
\dprod{k,i}{}\ 
\dfrac{[y_k+x_i+\alpha_i\delta]_{\nu_k}}
{[y_k+x_i]_{\nu_k}},
\end{eqs}
where $i,j\in\{1,\ldots,m\}$ and 
$k,l\in\{1,\ldots,n\}$.
\end{thm}

\par\medskip 
The identity \eqref{eq:xy} can be thought of 
as a transformation formula for certain multiple 
hypergeometric series. 
In order to clarify the idea, 
we introduce the multiple hypergeometric series 
$\Phi^{m,n}_N$ with respect to $[x]$ as follows. 
Given four vectors of variables
\begin{equation}
(a_1,\ldots,a_m),
(x_1,\ldots,x_m)\in\mathbb{C}^m,
\ \  \mbox{and}\ \ 
(b_1,\ldots,b_n),
(c_1,\ldots,c_n)\in\mathbb{C}^n,
\end{equation}
we define 
\begin{eqs}{ll}\label{eq:defPhi}
\medskip
\MP{m,n}{N}
{a_1,\ldots,a_m\\ x_1,\ldots,x_m}
{b_1,\ldots,b_n\\ c_1,\ldots,c_n}
\\
=\dsum{\mu\in\mathbb{N}^m,\,|\mu|=N}{}\ 
\dfrac{\Delta(x+\mu\delta)}{\Delta(x)}
\\
\qquad\quad\cdot
\dprod{1\le i,j\le m}\ 
\dfrac{[x_i-x_j+a_j]_{\mu_i}}{[x_i-x_j+\delta]_{\mu_i}}
\dprod{1\le i\le m; 1\le k\le n}\ 
\dfrac{[x_i+b_k]_{\mu_i}}{[x_i+c_k]_{\mu_i}}
\end{eqs}
for each $N=0,1,2,\ldots$.
Here the summation is taken over the finite set of 
$\mu\in\mathbb{N}^m$ with $|\mu|=N$.
Note that $\Phi^{m,n}_N$ defined as above 
is invariant with respect to the simultaneous permutations 
of $(a_1,\ldots,a_m)$ and $(x_1,\ldots,x_m)$, 
and symmetric in the variables $(b_1,\ldots,b_n)$, 
and in $(c_1,\ldots,c_n)$, respectively. 
We also remark that our parametrization of $\Phi^{m,n}_N$ is 
redundant in the sense that 
\begin{eqs}{ll}
\smallskip
\MP{m,n}{N}{a_1,&\ldots,&a_m \\ x_1+t,&\ldots,&x_m+t}
{b_1,\ldots,b_n \\ c_1,\ldots,c_n}
\\
=\MP{m,n}{N}{a_1,\ldots,a_m \\ x_1,\ldots,x_m}
{b_1+t,\ldots,b_n+t \\ c_1+t,\ldots,c_n+t}. 
\end{eqs}
We keep this redundancy for the sake of symmetry. 
\par\medskip
In this notation of $\Phi^{m,n}_N$, Theorem \ref{thm:PPd} can 
be understood as the identity 
\begin{eqs}{ll}
\medskip
\MP{m,n}{N}
{a_1,\ldots,a_m\\ x_1,\ldots,x_m}
{y_1-b_1,&\ldots,&y_n-b_n\\ y_1,&\ldots,&y_n}
\\
=
\MP{n,m}{N}
{b_1,\ldots,b_n\\ y_1,\ldots,y_n}
{x_1-a_1,&\ldots,&x_m-a_m\\x_1,&\ldots,&x_m}
\end{eqs}
for the special values 
\begin{eqs}{ll}
a_i=-\alpha_i\delta\ \ (i=1,\ldots,m) \ \ \mbox{and}\ \  
b_k=-\beta_k\delta\ \ (k=1,\ldots,n) 
\end{eqs}
for all 
$\alpha\in\mathbb{N}^m$, $\beta\in\mathbb{N}^n$ such that 
$|\alpha|=|\beta|.$
Theorem \ref{thm:PPd} is then generalized to 
the following transformation formula 
between $\Phi^{m,n}_N$ and $\Phi^{n,m}_N$. 

\begin{thm}[Duality transformation]
\label{thm:PPN}
Suppose that 
$(a_1,\ldots,a_m)$ $\in\mathbb{C}^m$ 
and $(b_1,\ldots,b_n)\in\mathbb{C}^n$ 
satisfy the balancing condition 
\begin{eqs}{ll}
a_1+\cdots+a_m=b_1+\cdots+b_n.
\end{eqs}
Then, for two sets of variables 
$(x_1,\ldots,x_m)$ and $(y_1,\ldots,y_n)$,
the following transformation formula holds 
between the multiple hypergeometric series 
$\Phi^{m,n}_N$ and $\Phi^{n,m}_N$\,$:$  
\begin{eqs}{ll}\label{eq:PPN}
\medskip
\MP{m,n}{N}
{a_1,&\ldots,&a_m\\ x_1,&\ldots,&x_m}
{y_1-b_1,&\ldots,&y_n-b_n\\ y_1,&\ldots,&y_n}
\\
=
\MP{n,m}{N}
{b_1,&\ldots,&b_n\\ y_1,&\ldots,&y_n}
{x_1-a_1,&\ldots,&x_m-a_m\\x_1,&\ldots,&x_m}
\end{eqs}
for each $N=0,1,2,\ldots$. 
\end{thm}
We remark that Theorem \ref{thm:PPN} is equivalent 
to the transformation formula
\begin{eqs}{ll}\label{eq:TFPP}
\smallskip
\MP{m,n}{N}
{a_1,\ldots,a_m\\ x_1,\ldots,x_m}
{b_1,\ldots,b_n\\ c_1,\ldots,c_n}
\\
=
\MP{n,m}{N}
{c_1-b_1,&\ldots,&c_n-b_n\\ c_1,&\ldots,&c_n}
{x_1-a_1,&\ldots,&x_m-a_m\\ x_1,&\ldots,&x_m}
\end{eqs}
for $N=0,1,2,\ldots$, 
under the balancing condition 
$a_1+\cdots+a_m+b_1+\cdots+b_n=c_1+\cdots+c_n$, 
or to 
\begin{eqs}{ll}\label{eq:bE}
\smallskip
\MP{m,n}{N}{a_1,&\ldots,&a_m \\ x_1,&\ldots,&x_m}
{b_1+y_1,&\ldots,&b_n+y_n\\ c+y_1,&\ldots,&c+y_n}
\\
=
\MP{n,m}{N}{c-b_1,&\ldots,&c-b_n\\ y_1,&\ldots,&y_n}
{c-a_1+x_1,&\ldots,&c-a_m+x_m\\ c+x_1,&\ldots,&c+x_m}
\end{eqs}
under the balancing condition 
$a_1+\cdots+a_m+b_1+\cdots+b_n=n c$.
Note that this transformation formula makes sense 
only in the setting of multiple hypergeometric series;
it becomes tautological when $(m,n)=(1,1)$. 
In what follows, 
we refer to the transformation formula \eqref{eq:PPN} 
as the {\em duality transformation} between 
$\Phi^{m,n}_N$ and $\Phi^{n,m}_N$.  
As we will see in the next section, 
this duality transformation formula implies 
various transformation formulas and summation 
formulas for very well-poised multiple elliptic 
hypergeometric series. 
\par\medskip
For the proof of Theorem \ref{thm:PPN}, 
we first remark that equality \eqref{eq:PPN} 
is reduced to the three typical cases, 
$[x]=x, \sin(\pi x)$ and $\sigma(x;\omega_1,\omega_2)$. 
In fact, if we replace $[x]$ 
by $e^{ax^2+b}[x]$, we obtain the same exponential factor
from the both sides of \eqref{eq:PPN}. 
Assuming that 
$[x]=x, \sin(\pi x)$ or $\sigma(x;\omega_1,\omega_2)$, 
we investigate the dependence of 
\begin{eqs}{ll}\smallskip
\mbox{LHS}=&\dsum{|\mu|=N}{}\ 
\dfrac{\Delta(x+\mu\delta)}{\Delta(x)}\, 
\dprod{i,j}\,
\dfrac{[x_i-x_j+a_j]_{\mu_i}}{[x_i-x_j+\delta]_{\mu_i}}
\dprod{i,k}\,
\dfrac{[x_i+y_k-b_k]_{\mu_i}}{[x_i+y_k]_{\mu_i}}
\end{eqs}
and 
\begin{eqs}{ll}\smallskip
\mbox{RHS}=&\dsum{|\nu|=N}{}\ 
\dfrac{\Delta(y+\nu\delta)}{\Delta(y)}\, 
\dprod{k,l}\,
\dfrac{[y_k-y_l+b_l]_{\nu_k}}{[y_k-y_l+\delta]_{\nu_k}}
\dprod{k,i}\,
\dfrac{[y_k+x_i-a_i]_{\nu_k}}{[y_k+x_i]_{\nu_k}}
\end{eqs}
on the variables $a_i$ and $b_k$.
In view of $a_1+\cdots+a_m=b_1+\cdots+b_n$, 
we consider the $(m+n-1)$-dimensional affine space 
$\mathbb{C}^{m+n-1}$ 
with coordinates 
$(a_1,\ldots,a_m,b_1,\ldots,b_{n-1})$, 
regarding $b_n$ as a linear function on this affine space. 
Let $F=F(a_1,\ldots,a_m,b_1,\ldots,b_{n-1})$ be 
the difference of the two sides regarded as a function in 
the variables $(a_1,\ldots,a_m,b_1,\ldots,b_{n-1})$. 
In the case where $[x]=x$,
$F$ is a polynomial in $a_i, b_k$. 
Also, in the case where $[x]=\sin(\pi x)$,
it is a Laurent polynomial in 
$e^{\pi\sqrt{-1} a_i},e^{\pi\sqrt{-1} b_k} $. 
We already know by 
Theorem \ref{thm:PPd} that 
$F$ has zeros at $a_i=-\alpha_i\delta$ 
($i=1,\ldots,m$), 
$b_k=-\beta_k\delta$ ($k=1,\ldots,n-1$) for all 
$(\alpha_1,\ldots,\alpha_m,\beta_1,\ldots,\beta_{n-1})
\in\mathbb{N}^{m+n-1}$. 
Under the assumption $[l\delta]\ne 0$ ($l\in\mathbb{Z}, l\ne 0$), 
we conclude by induction on the number of variables 
that $F$ is identically zero. 
Namely, the both sides coincide identically 
as holomorphic functions on $\mathbb{C}^{m+n-1}$
(and also as meromorphic functions in all the variables). 

In the case where $[x]$ is elliptic, 
Theorem \ref{thm:PPN} is proved by using the following 
lemma. 
\begin{lem}\label{lem:qper}
Let $F(z_1,\ldots,z_N)$ be a holomorphic function 
on $\mathbb{C}^N$.
Let $\omega_1,\omega_2\in\mathbb{C}$ be linearly 
independent over $\mathbb{R}$, and 
suppose that $F(z_1,\ldots,z_N)$ is quasi-periodic 
in each variable $z_i$ $(i=1,\ldots,N)$ with respect to the 
lattice $L=\mathbb{Z}\,\omega_1+\mathbb{Z}\,\omega_2$ 
in the following sense\,$:$ 
\begin{equation} 
F(z_1,\ldots,z_i+\omega_k,\ldots,z_N)
=e^{f_{i,k}(z_i)}\,F(z_1,\ldots,z_N)
\end{equation}
for some linear functions $f_{ik}(u)=a_{ik}u+b_{ik}$ 
$(i=1,\ldots,N; k=1,2)$. 
Let $\delta\in\mathbb{C}$ be a constant such that  
$\delta\not\in\mathbb{Q}\,\omega_1+\mathbb{Q}\,\omega_2$,
and suppose that 
\begin{equation}
F(\nu_1\delta,\ldots,\nu_N\delta)=0
\end{equation}
for all $\nu=(\nu_1,\ldots,\nu_N)\in\mathbb{N}^N$. 
Then $F(z_1,\ldots,z_N)$ is identically zero 
on $\mathbb{C}^N$. 
\end{lem}
This lemma is reduced by the induction on $N$ 
to the fact in the one variable case 
that any quasi-periodic holomorphic function 
may have only a finite number of zeros in the 
period-parallelogram. 

Let us regard the both sides of \eqref{eq:PPN} 
as holomorphic functions in the variables 
$(a_1,\ldots,a_m,b_1,\ldots,b_{n-1})\in\mathbb{C}^{m+n-1}$. 
Then one can check directly that they 
are quasi-periodic with the same multiplicative 
factors, with respect to each of the coordinate functions 
of $\mathbb{C}^{m+n-1}$. 
We know by Theorem \ref{thm:PPd} that 
the difference of the two sides has zeros 
at $a_i=-\alpha_i\delta$
($i=1,\ldots,m$), 
$b_j=-\beta_k\delta$ ($k=1,\ldots,n-1$)
for all $(\alpha_1,\ldots,\alpha_m,
\beta_1,\ldots,\beta_{n-1})\in\mathbb{N}^{m+n-1}$. 
This is impossible unless the both sides 
coincide identically as holomorphic functions
by Lemma \ref{lem:qper}.
This completes the proof of Theorem \ref{thm:PPN}. 

\begin{rem}\label{rem:phi}\rm\ \ 
In the trigonometric case, 
our multiple hypergeometric series $\Phi^{m,n}_N$ 
gives rise to the multiple $q$-hypergeometric series
\begin{eqs}{ll}
\smallskip
\Mphi{m,n}{N}
{a_1,\ldots,a_m\\ x_1,\ldots,x_m}
{b_1,\ldots,b_n\\ c_1,\ldots,c_n}
\\ \smallskip
=\dsum{\mu\in\mathbb{N}^m; |\mu|=N}{}\ 
\dprod{i<j}{}
\dfrac{q^{\mu_i}x_i-q^{\mu_j}x_j}{x_i-x_j}\ 
\dprod{i,j}{}
\dfrac{(a_jx_i/x_j;q)_{\mu_i}}
{(qx_i/x_j;q)_{\mu_i}}\ 
\dprod{i,k}{}
\dfrac{(b_kx_i;q)_{\mu_i}}
{(c_kx_i;q)_{\mu_i}},
\end{eqs}
where $i,j\in\{1,\ldots,m\}$ and $k\in\{1,\ldots,n\}$.
Setting $x=e^{2\pi\sqrt{-1}\xi}$, 
let us consider the case where 
\begin{eqs}{ll}
[\xi]=e^{\pi\sqrt{-1}\xi}-e^{-\pi\sqrt{-1}\xi}
=x^{\frac{1}{2}}-x^{-\frac{1}{2}}. 
\end{eqs}
Then the shifted factorials $[\xi]_k$ $(k=0,1,2,\ldots)$ are rewritten 
into 
\begin{eqs}{ll}
[\xi]_k=(-x^{-\frac{1}{2}})^k q^{-\frac{1}{2}{k \choose 2}}(x;q)_k,
\quad(x;q)_k=\dprod{i=0}{k-1}(1-xq^i)
\end{eqs}
with base $q=e^{2\pi\sqrt{-1}\delta}$. 
In this setting, $\Phi^{m,n}_{N}$ and $\phi^{m,n}_{N}$ 
are related by the formula
\begin{eqs}{ll}
\MP{m,n}{N}{\alpha_1,\ldots,\alpha_m \\ \xi_1,\ldots,\xi_m}
{\beta_1,\ldots,\beta_n\\ \gamma_1,\ldots,\gamma_n} \\
=
\left(\dfrac{qc_1\cdots c_n}
{a_1\cdots a_m b_1\cdots b_n}\right)^{\frac{N}{2}}
\Mphi{m,n}{N}
{a_1,\ldots,a_m \\ x_1,\ldots,x_m}
{b_1,\ldots,b_n\\ c_1,\ldots,c_n}, 
\end{eqs}
where 
\begin{eqs}{lll}
a_i=e^{2\pi\sqrt{-1}\alpha_i},\ \ 
&x_i=e^{2\pi\sqrt{-1}\xi_i}\ \ &(i=1,\ldots,m),\ \ \\
b_k=e^{2\pi\sqrt{-1}\beta_k},\ \ 
&c_k=e^{2\pi\sqrt{-1}\gamma_k}\ \ &(k=1,\ldots,n). 
\end{eqs}
Hence Theorem \ref{thm:PPN} implies the transformation formula
\begin{eqs}{ll}\label{eq:ppN}
\medskip
\Mphi{m,n}{N}
{a_1,&\ldots,&a_m\\ x_1,&\ldots,&x_m}
{y_1/b_1,&\ldots,&y_n/b_n\\ y_1,&\ldots,&y_n}
\\
=
\Mphi{n,m}{N}
{b_1,&\ldots,&b_n\\ y_1,&\ldots,&y_n}
{x_1/a_1,&\ldots,&x_m/a_m\\x_1,&\ldots,&x_m}
\end{eqs}
for each $N=0,1,2,\ldots$, 
under the balancing condition $a_1\cdots a_m=b_1\cdots b_n$. 
This formula can also be written in the form 
\begin{eqs}{l}
\medskip
\Mphi{m,n}{N}
{a_1,&\ldots,&a_m\\ x_1,&\ldots,&x_m}
{b_1 y_1,&\ldots,&b_n y_n\\  c y_1,&\ldots,&c y_n}
\\
=
\Mphi{n,m}{N}
{c/b_1,&\ldots,&c/b_n\\ y_1,&\ldots,&y_n}
{cx_1/a_1,&\ldots,& cx_m/a_m\\c x_1,&\ldots,&c x_m}
\end{eqs}
for $N=0,1,2,\ldots$ under the balancing condition 
$a_1\cdots a_m b_1\cdots b_n=c^n$. 

As is already known in \cite{K}, 
the duality transformation formula \eqref{eq:ppN}
for multiple basic hypergeometric series 
can be further generalized to non-balanced cases. 
Introducing a variable $u$, consider the infinite series
\begin{equation}
\Mphiu{m,n}
{a_1,&\ldots,&a_m\\ x_1,&\ldots,&x_m}
{b_1,&\ldots,&b_n\\ c_1,&\ldots,&c_n}
{u}
=\dsum{N=0}{\infty} u^N
\Mphi{m,n}{N}
{a_1,&\ldots,&a_m\\ x_1,&\ldots,&x_m}
{b_1,&\ldots,&b_n\\ c_1,&\ldots,&c_n}.
\end{equation}
Then the following {\em Euler transformation formula} 
holds for general parameters $a_i$, $b_k$ and $c$
(\cite{K}, Theorem 1.1):
\begin{eqs}{l}
\medskip
\Mphiu{m,n}
{a_1,&\ldots,&a_m\\ x_1,&\ldots,&x_m}
{b_1 y_1,&\ldots,&b_n y_n\\  c y_1,&\ldots,&c y_n}
{u}
\\
=
\dfrac{(a_1\cdots a_m b_1\cdots b_n u / c^n; q)_\infty}
{(u; q)_\infty}
\\
\quad
\Mphiu{n,m}
{c/b_1,&\ldots,&c/b_n\\ y_1,&\ldots,&y_n}
{cx_1/a_1,&\ldots,& cx_m/a_m\\c x_1,&\ldots,&c x_m}
{a_1\cdots a_m b_1\cdots b_n u/c^n}.
\end{eqs}
It would be an interesting problem to find an 
elliptic extension of this Euler transformation formula. 
\end{rem}

\section{Very well-poised multiple series $E^{m,n}$}
\setcounter{equation}{0}
In this section we show that our 
duality transformation between 
$\Phi^{m,n}_N$ and $\Phi^{n,m}_{N}$ 
(Theorem \ref{thm:PPN}) implies 
various transformation and summation formulas for 
{\em very well-poised} multiple elliptic hypergeometric 
series.

\par\medskip
We first remark that, when $m=1$, 
$\Phi^{m,n}_N$ reduces to a single term:
\begin{equation}\label{eq:Phi1n}
\renewcommand{\arraystretch}{0.8}\arraycolsep=2pt
\Phi^{1,n}_N\Big(
\begin{array}{c} a\\ x\end{array}\Big|
\begin{array}{c} b_1,\ldots,b_n\\ c_1,\ldots,c_n\end{array}
\Big)
=\dfrac{[a]_N}{[\delta]_N}\,\dfrac{[x+b_1]_N\cdots[x+b_n]_N}
{[x+c_1]_N\cdots[x+c_n]_N}. 
\end{equation}
In general, $\Phi^{1+m,n}_N$ is rewritten into 
a terminating $m$-dimensional sum which can be 
regarded as a {\em very well-poised} terminating 
multiple hypergeometric series:
\begin{eqs}{ll}\label{eq:Phi1+mn}
\smallskip
\MP{1+m,n}{N}
{a_0,a_1,\ldots,a_m\\ x_0,x_1,\ldots,x_m}
{b_1,\ldots,b_n\\ c_1,\ldots,c_n}
\\ \smallskip
=
\dfrac{[a_0]_N}{[\delta]_N}
\dprod{i=1}{m}\dfrac{[x_0-x_i+a_i]_N}{[x_0-x_i]_N}
\dprod{k=1}{n}\dfrac{[x_0+b_k]_N}{[x_0+c_k]_N}
\\ \smallskip
\quad\cdot
\dsum{\mu\in\mathbb{N}^m}{}\ 
\dfrac{\Delta(x_1+\mu_1\delta,\ldots,
x_m+\mu_m\delta)}{\Delta(x_1,\ldots,x_m)}
\,
\dprod{i=1}{m}
\dfrac{[x_i-x_0-N\delta+(|\mu|+\mu_i)\delta]}
{[x_i-x_0-N\delta]}
\\ \smallskip
\qquad\quad\cdot
\dfrac{[-N\delta]_{|\mu|}}
{[(1-N)\delta-a_0]_{|\mu|}}
\dprod{i=1}{m}
\dfrac{[x_i-x_0+a_0]_{\mu_i}}
{[x_i-x_0+\delta]_{\mu_i}}
\\ \smallskip
\qquad\quad\cdot
\dprod{j=1}{m}\ \left(
\dfrac{[-N\delta-x_0+x_j]_{|\mu|}}
{[(1-N)\delta-x_0+x_j-a_j]_{|\mu|}}\,
\dprod{i=1}{m}
\dfrac{[x_i-x_j+a_j]_{\mu_i}}
{[x_i-x_j+\delta]_{\mu_i}}\right)
\\ \smallskip
\qquad\quad\cdot
\dprod{k=1}{n}\ \left(
\dfrac{[(1-N)\delta-x_0-c_k]_{|\mu|}}
{[(1-N)\delta-x_0-b_k]_{|\mu|}}\,
\dprod{i=1}{m}
\dfrac{[x_i+b_k]_{\mu_i}}
{[x_i+c_k]_{\mu_i}}\right). 
\end{eqs}
This expression is obtained simply by applying the 
formula 
\begin{eqs}{ll}
[a]_{N-r}=\dfrac{[a]_N}{[a+(N-r)\delta]\cdots[a+(N-1)\delta]}
=\dfrac{(-1)^r [a]_N}{[(1-N)\delta-a]_r}
\end{eqs}
to every appearance of the shifted factorial of the form 
$[a]_{\mu_0}$, $\mu_0=N-(\mu_1+\cdots+\mu_n)$. 
Notice that by definition \eqref{eq:defPhi}
this series is symmetric with respect to the 
simultaneous permutations of the variables 
$(a_0,a_1,\ldots,a_m)$ and $(x_0,x_1,\ldots,x_m)$,
although it is not apparent in this expression. 

In view of the formula \eqref{eq:Phi1+mn}, 
we introduce the following 
notation for {\em very well-poised} 
multiple hypergeometric series $E^{m,n}$\,: 
\begin{eqs}{ll}\label{eq:defME}
\smallskip
\ME{m,n}{a_1,\ldots,a_m\\ x_1,\ldots,x_m}
{s}{u_1,\ldots,u_n}{v_1,\ldots,v_n}
\\ \smallskip
=
\dsum{\mu\in\mathbb{N}^m}{}\,
\dfrac{\Delta(x_1+\mu_1\delta,\ldots,x_m+\mu_m\delta)}
{\Delta(x_1,\ldots,x_m)}\,
\dprod{i=1}{m}
\dfrac{[x_i+s+(|\mu|+\mu_i)\delta]}{[x_i+s]}\,
\\ \smallskip
\qquad\cdot\ 
\dprod{j=1}{m}
\left(
\dfrac{[s+x_j]_{|\mu|}}{[\delta+s+x_j-a_j]_{|\mu|}}
\dprod{i=1}{m}
\dfrac{[x_i-x_j+a_j]_{\mu_i}}{[x_i-x_j+\delta]_{\mu_i}}
\right)
\\ \smallskip
\qquad\cdot\ 
\dprod{k=1}{n}
\left(
\dfrac{[v_k]_{|\mu|}}{[\delta+s-u_k]_{|\mu|}}
\dprod{i=1}{m}
\dfrac{[x_i+u_k]_{\mu_i}}{[x_i+\delta+s-v_k]_{\mu_i}}
\right).
\end{eqs}
This series contains 
``well-poised'' combinations of factors 
\begin{eqs}{ll}
\dfrac{[x_1+u]_{\mu_1}\cdots [x_m+u]_{\mu_m}}{[\delta+s-u]_{|\mu|}},
\ \ 
\dfrac{[v]_{|\mu|}}
{[\delta+s+x_1-v]_{\mu_1}\cdots[\delta+s+x_m-v]_{\mu_m}}
\end{eqs}
for $u=a_j-x_j, u_k$ and $v=s+x_j, v_k$. 
In order to make this series \eqref{eq:defME} terminate, 
we always work in such 
a situation where one of the following conditions is satisfied:
either 
\begin{equation}\label{eq:term}
\mbox{
\begin{tabular}{rl}\smallskip
(A)&\ $v_k=-N\delta$ \ for some $k=1,\ldots,n$ and 
$N=0,1,2,\ldots$; or \\
(B)&\ $a_i=-\alpha_i\delta$ \ ($i=1,\ldots,m$) \ 
for some $\alpha=(\alpha_1,\ldots,\alpha_m)\in\mathbb{N}^m$. 
\end{tabular}
}
\end{equation}
Note that our parametrization of $E^{m,n}$ has 
the following redundancy:
\begin{eqs}{ll}
\smallskip
\ME{m,n}{a_1,\ldots,a_m\\ x_1+t,\ldots,x_m+t}
{s}{u_1,\ldots,u_n}{v_1,\ldots,v_n}
\\ 
=
\ME{m,n}{a_1,\ldots,a_m\\ x_1,\ldots,x_m}
{s+t}{u_1+t,\ldots,u_n+t}{v_1,\ldots,v_n}.
\end{eqs}

When $m=1$, $E^{1,n}$ reduces to the very well-poised 
elliptic hypergeometric series ${}_{2n+4}E_{2n+3}$: 
\begin{eqs}{ll}
\smallskip
\ME{1,n}{a\\ x}{s}{u_1,\ldots,u_n}{v_1,\ldots,v_n}
\\
={}_{2n+4}E_{2n+3}\Big(x+s;a,x+u_1,\ldots,x+u_n,v_1,\ldots,v_n\Big),
\end{eqs}
where 
\begin{eqs}{ll}
{}_{r+1}E_{r}(s;u_1,\ldots,u_{r-2})
=\dsum{k=0}{\infty}
\dfrac{[s+2k\delta]}{[s]}
\dfrac{[s]_k}{[\delta]_k}\,
\dprod{i=1}{r-2}\dfrac{[u_i]_k}{[\delta+s-u_i]_k}. 
\end{eqs}
In the particular case where $m=1$ and $x=0$,
\begin{eqs}{ll}
\smallskip
\ME{1,n}{a\\ 0}{s}{u_1,\ldots,u_n}{v_1,\ldots,v_n}
\\
={}_{2n+4}E_{2n+3}\Big(s;a,u_1,\ldots,u_n,v_1,\ldots,v_n\Big). 
\end{eqs}
is completely symmetric in the $(2n+1)$ variables 
$(a,u_1,\ldots,u_n,v_1,\ldots,v_n)$.

It would be worthwhile to note here that 
$E^{m,n}$ carries a remarkable property 
concerning periodicity(cf.\,\cite{S}). 
We now use the abbreviated notation 
\begin{eqs}{ll}
\ME{m,n}{a \\ x}{s}{u}{v}
=\ME{m,n}{a_1,\ldots,a_m\\ x_1,\ldots,x_m}
{s}{u_1,\ldots,u_n}{v_1,\ldots,v_n}
\end{eqs}
for $a=(a_1,\ldots,a_m)$, $x=(x_1,\ldots,x_m)$,
$u=(u_1,\ldots,u_n)$, $v=(v_1,\ldots,v_n)$. 
\begin{lem}\label{lem:perE}
Suppose that the function $[x]$ is quasi-periodic with 
respect to a period $\omega\in\mathbb{C}$ in the sense that
\begin{eqs}{ll}
[x+\omega]\ =\ e^{\xi x+\eta} \ [x]\qquad(x\in \mathbb{C})
\end{eqs}
for some constants $\xi,\eta\in\mathbb{C}$. 
Then we have
\begin{eqs}{ll}
\ME{m,n}{a+l\,\omega \\ x}{s}{u+p\,\omega}{v+q\,\omega}
=
\ME{m,n}{a\\ x}{s}{u}{v}
\end{eqs}
for any $l\in\mathbb{Z}^m$, $p,q\in\mathbb{Z}^n$ such that 
$|l|+|p|+|q|=0$. 
\end{lem}

\par\medskip
\subsection*{\normalsize Relation between $\Phi^{m,n}_N$ and $E^{m,n}$}
\par\medskip
\begin{prop}\label{lem:PtoE}
The two types of multiple hypergeometric series 
$\Phi^{m, n}_N$ and $E^{m,n}$ are 
related as follows$:$ 
\begin{eqs}{ll}\label{eq:PtoE}
\smallskip
\MP{1+m,n}{N}
{a_0,a_1,\ldots,a_m\\ x_0,x_1,\ldots,x_m}
{b_1,\ldots,b_n\\ c_1,\ldots,c_n}
\\ \smallskip
=
\dfrac{[a_0]_N}{[\delta]_N}
\dprod{i=1}{m}\dfrac{[x_0-x_i+a_i]_N}{[x_0-x_i]_N}
\dprod{k=1}{n}\dfrac{[x_0+b_k]_N}{[x_0+c_k]_N}
\\ \smallskip
\quad\cdot\ 
E^{m,n+1}\Big(
\mat{a_1,\ldots,a_m\\x_1,\ldots,x_m}\,\Big|\,
{-N\delta-x_0}\,;{a_0-x_0,b_1,\ldots,b_n}\,;
\\
\qquad\qquad
{-N\delta, (1-N)\delta-x_0-c_1,\ldots,(1-N)\delta-x_0-c_n}\Big). 
\end{eqs}
Conversely, 
\begin{eqs}{ll}\label{eq:EtoP}
\smallskip
\ME{m,n+1}
{a_1,\ldots,a_m\\x_1,\ldots,x_m}
{s}{u_0,u_1,\ldots,u_n}{-N\delta,v_1,\ldots,v_n}
\\ \smallskip
=
\dfrac{[-N\delta]_N}{[\delta+s+u_0]_N}
\dprod{i=1}{m}\dfrac{[\delta+s+x_i]_N}{[\delta+s+x_i-a_i]_N}
\dprod{k=1}{n}\dfrac{[v_k]_N}{[\delta+s-u_k]_N}\,
\\ \smallskip
\quad\cdot\ 
\MP{1+m,n}{N}
{-N\delta-s+u_0,&a_1,\ldots,a_m\\-N\delta-s,&x_1,\ldots,x_m}
{u_1,&\ldots,&u_n\\ \delta+s-v_1,&\ldots,&\delta+s-v_n}.
\end{eqs}
\end{prop}
This proposition is an elliptic version of \cite{M1}, 
Lemma 1.22(see also \cite{M0}).
Note in particular that $\Phi^{2,n}_N$ reduces to a very well-poised 
terminating ${}_{2n+6}E_{2n+5}$:
\begin{eqs}{ll}\label{eq:Phi2n}
\smallskip
\MP{2,n}{N}
{a_0,a_1\\ x_0,x_1}
{b_1\ldots,b_n\\ c_1,\ldots,c_n}
\\ \smallskip
=
\dfrac{[a_0]_N[x_0-x_1+a_1]_N}{[\delta]_N[x_0-x_1]_N}
\dfrac{[x_0+b_1]_N\cdots[x_0+b_n]_N}
{[x_0+c_1]_N\cdots[x_0+c_n]_{N}}
\\ 
\quad\cdot\ 
{}_{2n+6}E_{2n+5}\big(
x_1-x_0-N\delta\ ;\ 
x_1-x_0+a_0, a_1, x_1+b_1,\ldots,x_1+b_n,
\\
\qquad\qquad\qquad\quad
-N\delta,(1-N)\delta-x_0-c_1,
\ldots,(1-N)\delta-x_0-c_n\big). 
\end{eqs}
Formula \eqref{eq:PtoE} is also valid for $m=0$ if 
we set $E^{0,n+1}=1$. 
\par\medskip
\subsection*{\normalsize Duality transformation between 
$E^{m,n+2}$ and $E^{n,m+2}$}
The duality transformation \eqref{eq:PPN} between $\Phi^{1+m,1+n}_N$ and 
$\Phi^{1+n,1+m}_N$ is now translated into a transformation 
formula between $E^{m,n+2}$ and $E^{n,m+2}$: 
\begin{eqs}{ll}\label{eq:EtoE}
\smallskip
\quad\dfrac{[\delta+s-c_1]_N[\delta+s-c_2]_N}{[d_1]_N[d_2]_N}
\dprod{i=1}{m}\dfrac{[\delta+s+x_i-a_i]_N}{[\delta+s+x_i]_N}
\dprod{k=1}{n}\dfrac{[\delta+s-u_k]_N}{[v_k]_N}\,
\\ \smallskip
\quad\cdot\ 
\ME{m,n+2}
{a_1,\ldots,a_m\\x_1,\ldots,x_m}
{s}{c_1,c_2,u_1,\ldots,u_n}{d_1,d_2,v_1,\ldots,v_n}
\\ \smallskip
=
\dfrac{[\delta+t+c_1]_N[\delta+t+c_2]_N}{[d_1]_N[d_2]_N}
\dprod{k=1}{n}\dfrac{[\delta+t+y_k-b_k]_N}{[\delta+t+y_k]_N}
\dprod{i=1}{m}\dfrac{[\delta+t-z_i]_N}{[w_i]_N}\,
\\ \smallskip
\quad\cdot\ 
\ME{n,m+2}
{b_1,\ldots,b_n\\y_1,\ldots,y_n}
{t}{-c_1,-c_2,z_1,\ldots,z_m}{d_1,d_2,w_1,\ldots,w_m}
\end{eqs}
under the balancing condition
\begin{eqs}{c}
\dsum{i=1}{m} a_i+\dsum{k=1}{2}(c_k+d_k)
+\dsum{k=1}{n}(u_k+v_k)
=(n+1)\delta+(n+2)s,
\end{eqs}
together with the termination condition 
$d_{2}=-N\delta$. 
Here the variables $b_k,y_k,t,z_i,w_i$ are specified by 
\begin{eqs}{lll}
t=d_1+d_2-s-\delta,\\
b_k=\delta+s-u_k-v_k,\ \ & y_k=\delta+s-v_k\quad&(k=1,\ldots,n),\\
z_i=x_i-a_i, &w_i=d_1+d_2-s-x_i\quad&(i=1,\ldots,m). 
\end{eqs}
Equivalently,
\begin{eqs}{ll}\label{eq:EtoE2}
\smallskip
\ME{m,n+2}
{a_1,\ldots,a_m\\x_1,\ldots,x_m}
{s}{c_1,c_2,u_1,\ldots,u_n}{d_1,d_2,v_1,\ldots,v_n}.
\\ \smallskip
=\dfrac{[\delta+s-c_1-d_1]_N[\delta+s-c_2-d_1]_N} 
{[\delta+s-c_1]_N[\delta+s-c_2]_N}
\\ \smallskip
\cdot\ 
\dprod{i=1}{m}\dfrac{[\delta+s+x_i]_N[\delta+s+x_i-a_i-d_1]_N}
{[\delta+s+x_i-a_i]_N[\delta+s+x_i-d_1]_N}
\,
\dprod{k=1}{n}\dfrac{[v_k]_N[\delta+s-u_k-d_1]_N}
{[\delta+s-u_k]_N[v_k-d_1]_N}
\,
\\ \smallskip
\cdot\ 
\ME{n,m+2}
{b_1,\ldots,b_n\\y_1,\ldots,y_n}
{t}{-c_1,-c_2,z_1,\ldots,z_m}{d_1,d_2,w_1,\ldots,w_m}. 
\end{eqs}
This transformation formula for the basic case 
was previously given in \cite{K}, and called  
the Bailey type transformation. 

\subsection*{\normalsize Dougall/Jackson summation formula for $E^{m,2}$}
When $n=0$, formula \eqref{eq:EtoE2} gives rise to 
the following summation formula
(\cite{R2}, Corollary 5.2): 
\begin{eqs}{ll}
\ME{m,2}
{a_1,\ldots,a_m\\x_1,\ldots,x_m}{s}{c_1,c_2}{d_1,d_2}
\\ \smallskip 
=
\dprod{k=1}{2}\dfrac{[\delta+s-c_k-d_1]_N}{[\delta+s-c_k]_N}
\dprod{i=1}{m}\dfrac{[\delta+s+x_i]_N[\delta+s+x_i-a_i-d_1]_N}
{[\delta+s+x_i-a_i]_N[\delta+s+x_i-d_1]_N},
\\ \smallskip
(a_1+\cdots+a_m+c_1+c_2+d_1+d_2=\delta+2s,\ \ d_2=-N\delta).
\end{eqs}
As the special case where $m=1$ and $x_1=0$,
this contains 
the Frenkel-Turaev summation formula 
(\cite{FT}, Theorem 5.5.2) for balanced ${}_{8}E_{7}$: 
\begin{eqs}{ll}
\smallskip
{}_{8}E_{7}\big(s;a,b,c,d,e\big)
\\ \smallskip
=
\dfrac{[\delta+s]_N[\delta+s-b-c]_N[\delta+s-b-d]_N
[\delta+s-c-d]_N}
{[\delta+s-b]_N[\delta+s-c]_N[\delta+s-d]_N[\delta+s-b-c-d]_N},
\\ \medskip
(a+b+c+d+e=\delta+2s,\quad e=-N\delta).
\end{eqs}

\subsection*{\normalsize 
Transformation from $E^{m,3}$ to ${}_{2m+8}E_{2m+7}$}

When $n=1$, the series $E^{1,m+2}$ 
in the right-hand side of \eqref{eq:EtoE2} 
gives a very well-poised elliptic hypergeometric series 
${}_{2m+8}E_{2m+7}$. 
Hence, our duality transformation 
provides a formula for rewriting balanced multiple elliptic 
hypergeometric series 
$E^{m,3}$ in terms of ${}_{2m+8}E_{2m+7}$: 
\begin{eqs}{ll}\label{eq:EtoE3}
\smallskip
\ME{m,3}
{a_1,\ldots,a_m\\x_1,\ldots,x_m}
{s}{c_0,c_1,c_2}{d_0,d_1,d_2}.
\\ \smallskip
=\dfrac{[d_0]_N[\delta+s-c_0-d_1]_N[\delta+s-c_1-d_1]_N[\delta+s-c_2-d_1]_N} 
{[d_0-d_1]_N[\delta+s-c_0]_N[\delta+s-c_1]_N[\delta+s-c_2]_N}
\\ \smallskip
\cdot\ 
\dprod{i=1}{m}\dfrac{[\delta+s+x_i]_N[\delta+s+x_i-a_i-d_1]_N}
{[\delta+s+x_i-a_i]_N[\delta+s+x_i-d_1]_N}
\\ \smallskip
\cdot\ 
{}_{2m+8}E_{2m+7}\big(t; d_1,d_2, e_0, e_1,e_2,
u_1,\ldots,u_m,v_1,\ldots,v_m\big)
\end{eqs}
under the condition 
\begin{eqs}{ll}
\dsum{i=1}{m} a_i+\dsum{k=0}{2}(c_k+d_k)=2\delta+3s,
\quad d_2=-N\delta, 
\end{eqs}
where 
\begin{eqs}{lll}
t=d_1+d_2-d_0,\quad &
e_k=\delta+s-d_0-c_k\ \ &(k=0,1,2),\\
u_i=\delta+s-d_0+x_i-a_i,\ \ &
v_i=d_1+d_2-s-x_i&(i=1,\ldots,m). 
\end{eqs}
This formula for the case where $m=1$ and $x_1=0$ recovers  
the following transformation formula for balanced ${}_{10}E_{9}$:
\begin{eqs}{ll}\label{eq:EtoE4}
\smallskip
{}_{10}E_{9}\big(
s; c_0,c_1,c_2,c_3,d_0,d_1,d_2\big)
\\ \smallskip
=
\dfrac{[d_0]_N[\delta+s]_N}{[d_0-d_1]_N[\delta+s-d_1]_N}
\dprod{k=0}{3}
\dfrac{[\delta+s-c_k-d_1]_N}{[\delta+s-c_k]_N}
\\ \smallskip
\cdot\ 
{}_{10}E_{9}\big(\widetilde{s};
\widetilde{c}_0,\widetilde{c}_1,
\widetilde{c}_2,
\widetilde{c}_3,\widetilde{d}_0,d_1,d_2\big)
\\ \smallskip
(c_0+c_1+c_2+c_3+d_0+d_1+d_2=2\delta+3s,\quad d_2=-N\delta), 
\end{eqs}
where 
\begin{eqs}{ll}
\widetilde{s}=d_1+d_2-d_0,\qquad
\widetilde{d}_0=d_1+d_2-s,\\
\widetilde{c}_k=\delta+s-d_0-c_k \quad(k=0,1,2,3).
\end{eqs}
 
\begin{rem}\label{rem:W}\rm\ \ 
In the trigonometric case where 
$[\xi]=x^{\frac{1}{2}}-x^{-\frac{1}{2}}$, $x=e^{2\pi\sqrt{-1}\xi}$,
the very well-poised multiple series $E^{m,n}$ gives 
\begin{eqs}{ll}
\smallskip
\MW{m,n}{a_1,\ldots,a_m\\ x_1,\ldots,x_m
}{s}{u_1,\ldots,u_n}{v_1,\ldots,v_n}
\\ \smallskip
=
\dsum{\mu\in \mathbb{N}^m}{}\ z^{|\mu|}
\dprod{1\le i<j\le m}{}\,
\dfrac{q^{\mu_i}x_i-q^{\mu_j}x_j}{x_i-x_j}\ 
\dprod{i=1}{m}
\dfrac{1-q^{|\mu|+\mu_i}sx_i}{1-sx_i}
\\ \smallskip
\quad\cdot\ 
\dprod{j=1}{m}
\dfrac{(sx_j;q)_{|\mu|}}{(qsx_j/a_j;q)_{|\mu|}}
\left(
\dprod{i=1}{m}
\dfrac{(a_jx_i/x_j;q)_{\mu_i}}{(qx_i/x_j;q)_{\mu_i}}
\right)
\\ 
\quad\cdot\ 
\dprod{k=1}{n}
\dfrac{(v_k;q)_{|\mu|}}{(qs/u_k;q)_{|\mu|}}
\left(
\dprod{i=1}{m}
\dfrac{(x_iu_k;q)_{\mu_i}}{(qsx_i/v_k;q)_{\mu_i}}
\right),
\end{eqs}
where $z=q^ns^n/a_1\cdots a_m u_1\cdots u_n v_1\cdots v_n$
(hence $z=q$ in the balanced case).
In fact we have
\begin{eqs}{ll}\smallskip
\ME{m,n}{\alpha_1,\ldots,\alpha_m\\ \xi_1,\ldots,\xi_m}
{\sigma}{\kappa_1,\ldots,\kappa_n}{\lambda_1,\ldots,\lambda_n}
\\
=
\MW{m,n}{a_1,\ldots,a_m\\ x_1,\ldots,x_m}
{s}{u_1,\ldots,u_n}{v_1,\ldots,v_n}
\end{eqs}
under the obvious correspondence of variables 
$a_i=e^{2\pi\sqrt{-1}\alpha_i},\ldots$.
\end{rem}

\section{Multiple Bailey transformations for $E^{m,3}$}
\setcounter{equation}{0} 
In this section, we derive two types of multiple generalizations 
of the elliptic Bailey transformation formulas 
due to Frenkel-Turaev \cite{FT}; one of them 
is the same as the one recently 
proved by Rosengren \cite{R2}, and the other seems to be new. 
\par\medskip
Recall that our duality transformation 
\eqref{eq:EtoE} transforms the balanced 
multiple series $E^{m,n+1}$'s 
into $E^{n,m+1}$. 
As we have seen in \eqref{eq:EtoE3}, 
in the particular case of $n=2$ 
it transforms the balanced $E^{m,3}$ into 
$E^{1,m+2}\propto{}_{2m+8}E_{2m+7}$; 
the corresponding  
${}_{2m+8}E_{2m+7}(s;u_1,\ldots,u_{2m+5})$ 
are completely symmetric with respect to 
the $2m+5$ parameters $u_1,\ldots,u_{2m+5}$. 
We can use this symmetry of ${}_{2m+8}E_{2m+7}$ 
to produce nontrivial transformation 
formulas for $E^{m,3}$ through the diagram 
\begin{eqs}{cccc}
E^{m,3} & \stackrel{\mbox{\small Bailey}}{\longrightarrow} & E^{m,3}\\
\downarrow & & \uparrow \\
{}_{2m+8}E_{2m+7} & \stackrel{\mbox{\small Symmetry}}
{\longrightarrow} &{}_{2m+8}E_{2m+7},
\end{eqs} 
where the vertical arrows represent the duality 
transformation \eqref{eq:EtoE3}. 
In this way, we obtain two types of 
Bailey transformation formulas for 
multiple elliptic hypergeometric series 
$E^{m,3}$. 
(A similar argument is used in \cite{K2} for 
studying multiple Sears' transformation formulas of type $A$ 
in the basic case.)

Let us write again the duality transformation formula 
between $E^{m,3}$ and ${}_{2m+8}E_{2m+7}$: 
Under the conditions 
\begin{eqs}{ll}\label{eq:bbtt}
\dsum{i=1}{m} a_i+\dsum{k=0}{2}(c_k+d_k)=2\delta+3s,\quad d_2=-N\delta, 
\end{eqs}
we have 
\begin{eqs}{ll}\label{eq:EtoE5}
\smallskip
\dfrac{[d_0-d_1]_N[\delta+s-c_0]_N[\delta+s-c_1]_N[\delta+s-c_2]_N}
{[d_0]_N[\delta+s-c_0-d_1]_N[\delta+s-c_1-d_1]_N[\delta+s-c_2-d_1]_N}\\ 
\smallskip
\cdot\ 
\dprod{i=1}{m}
\dfrac{[\delta+s+x_i-a_i]_N[\delta+s+x_i-d_1]_N}
{[\delta+s+x_i]_N[\delta+s+x_i-a_i-d_1]_N}
\\ \smallskip
\cdot\ 
\ME{m,3}
{a_1,\ldots,a_m\\x_1,\ldots,x_m}
{s}{c_0,c_1,c_2}{d_0,d_1,d_2}
\\ \smallskip
=
{}_{2m+8}E_{2m+7}\big(t; d_1,d_2, e_0, e_1,e_2,
u_1,\ldots,u_m,v_1,\ldots,v_m\big). 
\end{eqs} 
Here the variables $t$, $e_k$ $(k=0,1,2)$ and $u_i ,v_i$ ($i=1,\ldots,m$) 
are defined by 
\begin{eqs}{lll}
t=d_1+d_2-d_0,\quad &
e_k=\delta+s-d_0-c_k\ \ &(k=0,1,2),\\
u_i=\delta+s-d_0+x_i-a_i,\ \ &
v_i=d_1+d_2-s-x_i&(i=1,\ldots,m). 
\end{eqs}
Note that \eqref{eq:bbtt} implies 
\begin{eqs}{ll}
d_1+d_2+e_0+e_1+e_2+\dsum{i=1}{n}(u_i+v_i)=(m+1)\delta+(m+2)t. 
\end{eqs}

On the set of variables 
\begin{eqs}{lll}\label{eqs:vars}
(a_1,\ldots,a_m,x_1,\ldots,x_m,s,c_0,c_1,c_2,d_0,d_1,d_2), 
\end{eqs}
we first consider the following change of variables:
\begin{eqs}{lll}
\mbox{(I)}\quad& \widetilde{s}=\delta+2s-c_2-d_0-d_1,\quad
\widetilde{c}_2=\delta+s-d_0-d_1,\\
&\widetilde{d}_0=\delta+s-c_2-d_1,\qquad\quad\,
\widetilde{d}_1=\delta+s-c_2-d_0,\\
&\widetilde{c_0}=c_0,\quad 
\widetilde{c_1}=c_1,\quad \widetilde{d_2}=d_2,\\ 
& \widetilde{a}_i=a_i,\quad \widetilde{x_i}=x_i\quad(i=1,\ldots,m). 
\end{eqs}
For a given function $\varphi=\varphi(a_1,a_2,\ldots)$, 
we will denote by  
$\widetilde{\varphi}=\varphi(\widetilde{a}_1,\widetilde{a}_2,\ldots)$ 
the function obtained by replacing 
the variables $a_1,a_2,\ldots$ by 
$\widetilde{a}_1, \widetilde{a}_2,\ldots$. 
Then the variables for ${}_{2m+8}E_{2m+7}$ are transformed into 
\begin{eqs}{lll}
\widetilde{t}=t,\quad
&\widetilde{d}_1=e_2,\quad
&\widetilde{d}_2=d_2,\\
\widetilde{e}_0=e_0,\quad
&\widetilde{e}_1=e_1,\quad
&\widetilde{e}_2=d_1,\\
\widetilde{u}_i=u_i,\quad
&\widetilde{v}_i=v_i\quad&(i=1,\ldots,m), 
\end{eqs}
which is the transposition of $d_1$ and $e_2$. 
Hence the right-hand side of \eqref{eq:EtoE5} is 
invariant under this change of variables. 
In terms of $E^{m,3}$ on the left-hand side, 
this invariance implies the following multiple generalization 
of the Bailey transformation due to Rosengren 
(\cite{R2}, Corollary 8.2). 

\begin{thm}[Multiple Bailey transformation I] 
Under the balancing condition 
\begin{eqs}{ll}
\dsum{i=1}{m} a_i+\dsum{k=0}{2}(c_k+d_k)=2\delta+3s,
\end{eqs} 
the following identify holds for $d_2=-N\delta$ $(N=0,1,2,\ldots,)$\,$:$
\begin{eqs}{ll}\label{eq:BaileyI}
\smallskip
\ME{m,3}
{a_1,\ldots,a_m\\x_1,\ldots,x_m}
{\widetilde{s}}{c_0,c_1,\widetilde{c}_2}
{\widetilde{d}_0,\widetilde{d}_1,d_2}
\\ \smallskip
=
\dfrac{[\delta+s-c_0]_N[\delta+s-c_1]_N}
{[\delta+\widetilde{s}-c_0]_N[\delta+\widetilde{s}-c_1]_N}
\dprod{i=1}{m}
\dfrac{[\delta+s+x_i-a_i]_N[\delta+\widetilde{s}+x_i]_N}
{[\delta+s+x_i]_N[\delta+\widetilde{s}+x_i-a_i]_N}
\\ \smallskip
\cdot\ 
\ME{m,3}
{a_1,\ldots,a_m\\x_1,\ldots,x_m}
{s}{c_0,c_1,c_2}{d_0,d_1,d_2},
\end{eqs}
where 
\begin{eqs}{ll}
\widetilde{s}=\delta+2s-c_2-d_0-d_1,\quad
\widetilde{c}_2=\delta+s-d_0-d_1,\\
\widetilde{d}_0=\delta+s-c_2-d_1,\qquad\quad\,
\widetilde{d}_1=\delta+s-c_2-d_0. 
\end{eqs}
\end{thm}
Formula \eqref{eq:BaileyI}
can be written alternatively 
in the form 
\begin{eqs}{ll}\label{eq:BaileyIa}
\smallskip
\ME{m,3}
{a_1,\ldots,a_m\\x_1,\ldots,x_m}
{\widetilde{s}}{c_0,c_1,\widetilde{c}_2}
{\widetilde{d}_0,\widetilde{d}_1,d_2}
\\ \smallskip
=
\dfrac{[\delta+s-c_0]_N[\delta+s-c_1]_N}
{[\delta+s-c_0-|a|]_N
[\delta+s-c_1-|a|]_N}
\\ \smallskip
\quad\cdot\ 
\dprod{i=1}{m}
\dfrac{[\delta+s+x_i-a_i]_N[\delta+s-x_i-c_0-c_1-|a|]_N}
{[\delta+s+x_i]_N[\delta+s-x_i+a_i-c_0-c_1-|a|]_N}
\\ \smallskip
\quad\cdot\ 
\ME{m,3}
{a_1,\ldots,a_m\\x_1,\ldots,x_m}
{s}{c_0,c_1,c_2}{d_0,d_1,d_2},
\end{eqs}
where $|a|=a_1+\cdots+a_m$. 

\par\medskip
Next we consider the following change of variables:
\begin{eqs}{lll}
\mbox{(II)}\quad&
\widetilde{s}=\delta+2s-c_0-c_1-c_2, \quad& 
\widetilde{c}_0=\delta+s-c_1-c_2,\\
&
\widetilde{c}_1=\delta+s-c_0-c_2,\quad&
\widetilde{c}_2=\delta+s-c_0-c_1,\\
&\widetilde{d}_k=d_k\quad(k=0,1,2),\\
& 
\widetilde{a}_i=a_i,\quad
\widetilde{x}_i=a_i-x_i-|a|\quad&(i=1,\ldots,m),
\end{eqs}
where $|a|=a_1+\cdots+a_m$. 
Then the variables for ${}_{2m+8}E_{2m+7}$ are 
transformed into 
\begin{eqs}{lll}
\widetilde{t}=t,\quad
&\widetilde{d}_1=d_1,\quad
&\widetilde{d}_2=d_2,\\
\widetilde{e}_k=e_k\quad
&(k=0,1,2),\\
\widetilde{u}_i=v_i,\quad
&\widetilde{v}_i=u_i\quad&(i=1,\ldots,m)
\end{eqs}
under the balancing condition \eqref{eq:bbtt}, 
which is the simultaneous transposition 
of $u_i$ and $v_i$ ($i=1,\ldots,m$). 
Since the right-hand side of \eqref{eq:EtoE5} is invariant 
under this transformation of variables, 
we obtain another transformation formula 
for $E^{m,3}$. 

\begin{thm}[Multiple Bailey transformation II] 
Under the balancing condition 
\begin{eqs}{ll}
\dsum{i=1}{m} a_i+\dsum{k=0}{2}(c_k+d_k)=2\delta+3s,
\end{eqs}
the following identity holds for $d_2=-N\delta$ $(N=0,1,2,\ldots)$\,$:$ 
\begin{eqs}{ll}\label{eq:BaileyII}
\smallskip
\ME{m,3}
{a_1,\ldots,a_m\\ \widetilde{x}_1,\ldots,\widetilde{x}_m}
{\widetilde{s}}{\widetilde{c}_0,\widetilde{c}_1,\widetilde{c}_2}
{d_0,d_1,d_2}
\\ \smallskip
=
\dprod{i=1}{m}
\dfrac{[\delta+s+x_i-a_i]_N[\delta+s+x_i-d_1]_N}
{[\delta+s+x_i]_N[\delta+s+x_i-a_i-d_1]_N}
\\ \smallskip
\cdot\ 
\dprod{i=1}{m}
\dfrac
{[\delta+\widetilde{s}+\widetilde{x}_i]_N
[\delta+\widetilde{s}+\widetilde{x}_i-a_i-d_1]_N}
{[\delta+\widetilde{s}+\widetilde{x}_i-a_i]_N
[\delta+\widetilde{s}+\widetilde{x}_i-d_1]_N}
\\ 
\cdot\ 
\ME{m,3}
{a_1,\ldots,a_m\\x_1,\ldots,x_m}
{s}{c_0,c_1,c_2}{d_0,d_1,d_2},
\end{eqs}
where 
\begin{eqs}{ll}
\widetilde{s}=\delta+2s-c_0-c_1-c_2, \quad& 
\widetilde{c}_0=\delta+s-c_1-c_2,\\
\widetilde{c}_1=\delta+s-c_0-c_2,\quad&
\widetilde{c}_2=\delta+s-c_0-c_1,\\
\widetilde{x}_i=a_i-x_i-|a|\quad&(i=1,\ldots,m).
\end{eqs}
\end{thm}
The second multiple Bailey transformation \eqref{eq:BaileyII} 
can also be written in the form 
\begin{eqs}{ll}\label{eq:BaileyIIa}
\smallskip
\ME{m,3}
{a_1,\ldots,a_m\\ \widetilde{x}_1,\ldots,\widetilde{x}_m}
{\widetilde{s}}{\widetilde{c}_0,\widetilde{c}_1,\widetilde{c}_2}
{d_0,d_1,d_2}
\\ \smallskip
=
\dprod{i=1}{m}
\dfrac{[\delta+s+x_i-d_0]_N[\delta+s+x_i-d_1]_N}
{[\delta+s+x_i]_N[\delta+s+x_i-d_0-d_1]_N}
\\ \smallskip
\cdot\ 
\dprod{i=1}{m}
\dfrac{[\delta+s+x_i-a_i]_N[\delta+s+x_i-a_i-d_0-d_1]_N}
{[\delta+s+x_i-a_i-d_0]_N[\delta+s+x_i-a_i-d_1]_N}
\\ 
\cdot\ 
\ME{m,3}
{a_1,\ldots,a_m\\x_1,\ldots,x_m}
{s}{c_0,c_1,c_2}{d_0,d_1,d_2}. 
\end{eqs}
We remark that, when $m=1$ (and $x_1=0$), both 
\eqref{eq:BaileyI} and \eqref{eq:BaileyII} reduce 
to the elliptic Bailey transformation formula 
due to Frenkel-Turaev (\cite{FT}, Theorem 5.5.1): 
\begin{eqs}{ll}
\smallskip
{}_{10}E_{9}\big(
\widetilde{s}; \widetilde{c}_0,\widetilde{c}_1,\widetilde{c}_2,
d_0,d_1,d_2,d_3
\big)
\\ \smallskip
=
\dfrac{[\delta+s-d_0]_N[\delta+s-d_1]_N
[\delta+s-d_2]_N[\delta+s-d_0-d_1-d_2]_N}
{[\delta+s]_N[\delta+s-d_0-d_1]_N
[\delta+s-d_0-d_2]_N[\delta+s-d_1-d_2]_N}
\\ \smallskip
\cdot\ 
{}_{10}E_{9}\big(
s; c_0,c_1,c_2,
d_0,d_1,d_2,d_3
\big),\\
\widetilde{s}=\delta+2s-c_0-c_1-c_2,\quad
\widetilde{c}_0=\delta+s-c_1-c_2,\\
\widetilde{c}_1=\delta+s-c_0-c_2,\qquad\quad\,
\widetilde{c}_2=\delta+s-c_0-c_1,
\end{eqs}
where 
$c_0+c_1+c_2+d_0+d_1+d_2+d_3=2\delta+3s,\quad d_3=-N\delta$. 

\par\bigskip 
Until now, we have discussed terminating multiple hypergeometric series 
\begin{eqs}{ll}
\ME{m,n}{a_1,\ldots,a_m\\ x_1,\ldots,x_m}
{s}{u_1,\ldots,u_n}{v_1,\ldots,v_n}
\end{eqs}
such that 
\begin{eqs}{ll}
\mbox{(A)}\quad & v_k=-N\delta\quad 
\mbox{for some}\quad
k=1,\ldots,n\ \ \mbox{and}\ \ N=0,1,2,\ldots.
\end{eqs}
As we remarked before, we can also consider terminating 
series under the condition 
\begin{eqs}{ll}
\mbox{(B)}\quad &a_i=-\alpha_i\delta\ \ 
(i=1,\ldots,m)\ \ 
\mbox{for some}\ \ 
\alpha=(\alpha_1,\ldots,\alpha_m)\in \mathbb{N}^m. 
\end{eqs}
There is a standard procedure for deriving identities 
for terminating series of type (B) from 
those of type (A), and {\em vice versa}.  
In the following, 
we derive two types of multiple Bailey transformation 
formulas for terminating series of type (B). 
In the first multiple Bailey transformation \eqref{eq:BaileyIa}, 
besides the balancing condition
\begin{eqs}{l}
\dsum{i=1}{m} a_i+\dsum{k=0}{2}(c_k+d_k)=2\delta+3s 
\end{eqs}
we assume that both the two termination 
conditions 
\begin{eqs}{ll}
a_i=-\alpha_i\delta\quad(i=1,\ldots,m)\quad
\mbox{and}\quad 
d_2=-N\delta
\end{eqs}
are satisfied. 
Then by using the identity 
\begin{eqs}{l}
\dfrac{[x+k\delta]_l}{[x]_l}=\dfrac{[x+l\delta]_k}{[x]_k}
\qquad(k,l\in\mathbb{N}), 
\end{eqs}
we can rewrite \eqref{eq:BaileyIa} into 
\begin{eqs}{ll}\label{eq:BaileyIb}
\smallskip
\ME{m,3}
{a_1,\ldots,a_m\\x_1,\ldots,x_m}
{\widetilde{s}}{c_0,c_1,\widetilde{c}_2}
{\widetilde{d}_0,\widetilde{d}_1,d_2}
\\ \smallskip
=
\dfrac{[\delta+s-c_0]_{|\alpha|}[\delta+s-c_1]_{|\alpha|}}
{[\delta+s-c_0-d_2]_{|\alpha|}
[\delta+s-c_1-d_2]_{|\alpha|}}
\\ \smallskip
\quad\cdot\ 
\dprod{i=1}{m}
\dfrac{[\delta+s+x_i-d_2]_{\alpha_i}
[\delta+s-x_i+a_i-|a|-c_0-c_1-d_2]_{\alpha_i}}
{[\delta+s+x_i]_{\alpha_i}[\delta+s-x_i+a_i-|a|-c_0-c_1]_{\alpha_i}}
\\ \smallskip
\quad\cdot\ 
\ME{m,3}
{a_1,\ldots,a_m\\x_1,\ldots,x_m}
{s}{c_0,c_1,c_2}{d_0,d_1,d_2},
\end{eqs}
where $|a|=a_1+\cdots+a_m$. 
This formula is valid for the values $d_2=-N\delta$ 
($N=0,1,2,\ldots$). 
We now regard each of the two-sides of this formula 
as a function of $d_2$; 
we regard $c_0$ as a linear function 
of $d_2$ in the form $c_0=\lambda-d_2$, 
keeping the other variables 
$a_i$, $x_i$ $(i=1,\ldots,m)$ and 
$s, c_1,c_2,d_0,d_1$ constant. 
In the elliptic case, it can be checked 
that the both sides of \eqref{eq:BaileyIb} 
are periodic in the variable $d_2$ in this sense
(see Lemma \ref{lem:perE}). 
Hence we conclude that \eqref{eq:BaileyIb} holds identically 
as a meromorphic function of $d_2$,  
only on the condition (B). 
Similarly, 
we can rewrite \eqref{eq:BaileyIIa} into 
\begin{eqs}{ll}\label{eq:BaileyIIb}
\smallskip
\ME{m,3}
{a_1,\ldots,a_m\\ \widetilde{x}_1,\ldots,\widetilde{x}_m}
{\widetilde{s}}{\widetilde{c}_0,\widetilde{c}_1,\widetilde{c}_2}
{d_0,d_1,d_2}
\\ \smallskip
=
\dprod{i=1}{m}
\dfrac{[\delta+s+x_i-d_0]_{\alpha_i}[\delta+s+x_i-d_1]_{\alpha_i}}
{[\delta+s+x_i]_{\alpha_i}[\delta+s+x_i-d_0-d_1]_{\alpha_i}}
\\ \smallskip
\cdot\ 
\dprod{i=1}{m}
\dfrac{[\delta+s+x_i-d_2]_{\alpha_i}[\delta+s+x_i-d_0-d_1-d_2]_{\alpha_i}}
{[\delta+s+x_i-d_0-d_2]_{\alpha_i}[\delta+s+x_i-d_1-d_2]_{\alpha_i}}
\\ 
\cdot\ 
\ME{m,3}
{a_1,\ldots,a_m\\x_1,\ldots,x_m}
{s}{c_0,c_1,c_2}{d_0,d_1,d_2}. 
\end{eqs}
By a similar reasoning for the variables $d_2$ and $d_0$, 
we see that this formula is valid under the balancing condition 
together with the termination condition (B). 
After all, we also obtain two types of multiple Bailey transformations 
under the condition (B). 

\begin{thm}
Suppose the parameters $a_i$ $(i=1,\ldots,m)$ and 
$c_k,d_k$ $(k=0,1,2)$ satisfy the balancing condition
\begin{eqs}{l}
\dsum{i=1}{m} a_i+\dsum{k=0}{2}(c_k+d_k)=2\delta+3s
\end{eqs}
and the termination condition 
\begin{eqs}{l}
a_i=-\alpha_i\delta\ \ (i=1,\ldots,m)\ \ \mbox{for some}\ \ 
\alpha=(\alpha_1,\ldots,\alpha_m)\in \mathbb{N}^m. 
\end{eqs}
Then we have the following two types of multiple Bailey 
transformations$:$ 
\begin{eqs}{ll}\label{eq:BaileyIaa}
\smallskip
\mbox{\rm(I)}\quad&\ME{m,3}
{a_1,\ldots,a_m\\x_1,\ldots,x_m}
{\widetilde{s}}{c_0,c_1,\widetilde{c}_2}
{\widetilde{d}_0,\widetilde{d}_1,d_2}
\\ \smallskip
&=
\dfrac{[\delta+s-c_0]_{|\alpha|}[\delta+s-c_1]_{|\alpha|}}
{[\delta+s-c_0-d_2]_{|\alpha|}
[\delta+s-c_1-d_2]_{|\alpha|}}
\\ \smallskip
&\quad\cdot\ 
\dprod{i=1}{m}
\dfrac{[\delta+s+x_i-d_2]_{\alpha_i}
[\delta+s-x_i+a_i-|a|-c_0-c_1-d_2]_{\alpha_i}}
{[\delta+s+x_i]_{\alpha_i}[\delta+s-x_i+a_i-|a|-c_0-c_1]_{\alpha_i}}
\\ \smallskip
&\quad\cdot\ 
\ME{m,3}
{a_1,\ldots,a_m\\x_1,\ldots,x_m}
{s}{c_0,c_1,c_2}{d_0,d_1,d_2},
\end{eqs}
where 
\begin{eqs}{ll}
\widetilde{s}=\delta+2s-c_2-d_0-d_1,\quad
\widetilde{c}_2=\delta+s-d_0-d_1,\\
\widetilde{d}_0=\delta+s-c_2-d_1,\qquad\quad\,
\widetilde{d}_1=\delta+s-c_2-d_0,
\end{eqs}
and
\begin{eqs}{ll}\label{eq:BaileyIIB}
\smallskip
\mbox{\rm(II)}\qquad&
\ME{m,3}
{a_1,\ldots,a_m\\ \widetilde{x}_1,\ldots,\widetilde{x}_m}
{\widetilde{s}}{\widetilde{c}_0,\widetilde{c}_1,\widetilde{c}_2}
{d_0,d_1,d_2}
\\ \smallskip
&=
\dprod{i=1}{m}
\dfrac{[\delta+s+x_i-d_0]_{\alpha_i}[\delta+s+x_i-d_1]_{\alpha_i}}
{[\delta+s+x_i]_{\alpha_i}[\delta+s+x_i-d_0-d_1]_{\alpha_i}}
\\ \smallskip
&\cdot\ 
\dprod{i=1}{m}
\dfrac{[\delta+s+x_i-d_2]_{\alpha_i}[\delta+s+x_i-d_0-d_1-d_2]_{\alpha_i}}
{[\delta+s+x_i-d_0-d_2]_{\alpha_i}[\delta+s+x_i-d_1-d_2]_{\alpha_i}}
\\ 
&\cdot\ 
\ME{m,3}
{a_1,\ldots,a_m\\x_1,\ldots,x_m}
{s}{c_0,c_1,c_2}{d_0,d_1,d_2},
\end{eqs}
where
\begin{eqs}{ll}
\widetilde{s}=\delta+2s-c_0-c_1-c_2, \quad& 
\widetilde{c}_0=\delta+s-c_1-c_2,\\
\widetilde{c}_1=\delta+s-c_0-c_2,\quad&
\widetilde{c}_2=\delta+s-c_0-c_1,\\
\widetilde{x}_i=a_i-x_i-|a|\quad&(i=1,\ldots,m). 
\end{eqs}
\end{thm}
(Formula \eqref{eq:BaileyIaa} is proved in \cite{R2}, Corollary 8.1.)
\begin{rem}\rm 
The second Bailey transformation formulas \eqref{eq:BaileyIIa} 
and \eqref{eq:BaileyIIB} appear to be new even in the 
case of multiple basic and ordinary hypergeometric series. 
In the basic case, 
under the balancing condition 
\begin{equation}
a_1\cdots a_m c_0c_1c_2 d_0d_1d_2=q^2s^3, 
\end{equation}
the change of parameters 
\begin{eqs}{ll}
\widetilde{s}=qs^2/c_0c_1c_2, \quad& 
\widetilde{c}_0=qs/c_1c_2,\\
\widetilde{c}_1=qs/c_0c_2,\quad&
\widetilde{c}_2=qs/c_0c_1,\\
\widetilde{x}_i=a_i/a_1\cdots a_m x_i\quad&(i=1,\ldots,m)
\end{eqs}
implies the following two transformation formulas 
for $W^{m,3}$ of Remark \ref{rem:W}: 
\begin{eqs}{l}
\smallskip
\MW{m,3}
{a_1,\ldots,a_m\\ \widetilde{x}_1,\ldots,\widetilde{x}_m}
{\widetilde{s}}{\widetilde{c}_0,\widetilde{c}_1,\widetilde{c}_2}
{d_0,d_1,d_2}
\\ \smallskip
=
\dprod{i=1}{m}
\dfrac{(qsx_i/d_0;q)_N(qsx_i/d_1;q)_N}
{(qsx_i;q)_N(qsx_i/d_0d_1;q)_N}
\\ \smallskip
\cdot\ 
\dprod{i=1}{m}
\dfrac{(qsx_i/a_i;q)_N(qsx_i/a_id_0d_1;q)_N}
{(qsx_i/a_id_0;q)_N(qsx_i/a_id_1;q)_N}
\\ 
\cdot\ 
\MW{m,3}
{a_1,\ldots,a_m\\x_1,\ldots,x_m}
{s}{c_0,c_1,c_2}{d_0,d_1,d_2},
\end{eqs}
for $d_2=q^{-N}$, $N=0,1,2\ldots$; and 
\begin{eqs}{l}
\smallskip
\MW{m,3}
{a_1,\ldots,a_m\\ \widetilde{x}_1,\ldots,\widetilde{x}_m}
{\widetilde{s}}{\widetilde{c}_0,\widetilde{c}_1,\widetilde{c}_2}
{d_0,d_1,d_2}
\\ \smallskip
=
\dprod{i=1}{m}
\dfrac{(qsx_i/d_0;q)_{\alpha_i}(qsx_i/d_1;q)_{\alpha_i}}
{(qsx_i;q)_{\alpha_i}(qsx_i/d_0d_1;q)_{\alpha_i}}
\\ \smallskip
\cdot\ 
\dprod{i=1}{m}
\dfrac{(qsx_i/d_2;q)_{\alpha_i}(qsx_i/d_0d_1d_2;q)_{\alpha_i}}
{(qsx_i/d_0d_2;q)_{\alpha_i}(qsx_i/d_1d_2;q)_{\alpha_i}}
\\ 
\cdot\ 
\MW{m,3}
{a_1,\ldots,a_m\\x_1,\ldots,x_m}
{s}{c_0,c_1,c_2}{d_0,d_1,d_2},
\end{eqs}
for 
$a_i=q^{-\alpha_i}$, $\alpha_i\in\mathbb{N}$ ($i=1,\ldots,m$).
\end{rem}

\end{document}